%% file: m3-II-4.tex
\input m3-macs

\pageno=239

\tinfo II.4.239-253

\SetTFLinebox{\gtp }
\SetFLinebox{\gtv3 }
\SetHLinebox{\issn}

\H 4. Drinfeld modules and\\
local fields of positive characteristic

 Ernst--Ulrich Gekeler

\SetAuthorHead{E.-U. Gekeler}
\SetTitleHead{Part II. Section 4. Drinfeld modules and 
local fields of positive characteristic \qquad\qquad}

The relationship between local fields and Drinfeld modules is twofold.
Drinfeld modules allow explicit construction of abelian and nonabelian 
extensions with prescribed properties of local and global fields
of positive characteristic. On the other hand, $n$-dimensional local
fields arise in the construction of (the compactification of) moduli
schemes $X$ for Drinfeld modules, such schemes being provided with a
natural stratification $X_0 \subset X_1 \subset \cdots X_i \cdots \subset
X_n = X$ through smooth subvarieties $X_i$ of dimen\-sion~$i$.

We will survey that correspondence, but refer to the literature for
detailed proofs (provided these exist so far). An important remark is
in order: The contents of this article take place in characteristic
$p > 0$, and are in fact locked up in the characteris\-tic~$p$ world. No
lift to characteristic zero nor even to schemes over $\Bbb Z/p^2$ 
is known!

\HH 4.1. Drinfeld modules

 Let $L$ be a field of characteristic $p$
containing the field $\Bbb F_q$, and denote by $\tau = \tau_q$
 raising to the $q$th power map $x \mapsto x^q$. If ``$a$''
denotes multiplication by $a \in L$, then $\tau a = a^q\tau$. The ring
${\rm End}({\Bbb G_a}_{/L})$ of endomorphisms of the additive group
${\Bbb G_a}_{/L}$ equals $L\{\tau_p\} = \{\sum a_i\tau_p^i:a_i\in L\}$, the
non-commutative polynomial ring in $\tau_p = (x \mapsto x^p)$ with the
above commutation rule $\tau_pa = a^p\tau$. Similarly, the subring
${\rm End}_{\Bbb F_q}({\Bbb G_a}_{/L})$ of $\Bbb F_q$-endomorphisms is
$L\{\tau\}$ with $\tau = \tau_p^n$ if $q = p^n$. Note that $L\{\tau\}$
is an $\Bbb F_q$-algebra since $\Bbb F_q \hookrightarrow L\{\tau\}$ is 
central.

\df Definition 1

 Let ${\Cal C}$ be a smooth geometrically connected
projective curve over $\Bbb F_q$.
Fix a closed (but not necessarily
$\Bbb F_q$-rational) point $\infty$ of ${\Cal C}$. The ring $A = 
\Gamma({\Cal C}-\{\infty\},{\Cal O}_{\Cal C})$ is called a {\it Drinfeld
ring}. Note that $A^{\ast} = \Bbb F_q^{\,\ast}$.
\enddf

\eg Example 1

 If ${\Cal C}$ is the projective line ${\Bbb P^1}_{/\Bbb F_q}$ 
and $\infty$ is the usual point at infinity then $A = \Bbb F_q[T]$.
\endeg
\eg Example 2

 Suppose that $p \not= 2$, that ${\Cal C}$ is given by an 
affine equation $Y^2 = f(X)$ with a separable polynomial $f(X)$ of
even positive degree with leading coefficient a non-square in $\Bbb F_q$, and that
$\infty$ is the point above $X = \infty$. Then $A = \Bbb F_q[X,Y]$ is
a Drinfeld ring with $\deg_{\Bbb F_q}(\infty)=2$.
\endeg

\df Definition 2

 An {\it $A$-structure} on a field $L$ is a homomorphism
of $\Bbb F_q$-algebras (in brief: an $\Bbb F_q$-ring homomorphism)
$\gamma\colon A \rightarrow L$. Its {\it $A$-characteristic} 
${\rm char}_A(L)$ is the maximal ideal ${\rm ker}(\gamma)$, if $\gamma$
fails to be injective, and $\infty$ otherwise. A {\it Drinfeld
module structure} on such a field $L$ is given by an $\Bbb F_q$-ring
homomorphism $\phi\colon  A \rightarrow L\{\tau\}$ such that
$\partial \circ \phi = \gamma$, where $\partial\colon  L\{\tau\} 
\rightarrow L$ is the $L$-homomorphism sending $\tau$ to 0.  

Denote $\phi(a)$ by $\phi_a \in {\rm End}_{\Bbb F_q}({\Bbb G_a}_{/L})$; 
$\phi_a$ induces on the additive group over $L$ (and on each $L$-algebra $M$) 
a {\it new} structure as an $A$-module: 
 $$a \ast x := \phi_a(x) \quad (a\in A, x \in M).\leqno{(4.1.1)}$$
We briefly call $\phi$ a Drinfeld module over $L$, usually omitting
reference to $A$.
\enddf

\df Definition 3

 Let $\phi$ and $\psi$ be Drinfeld modules over the 
$A$-field $L$. A {\it homomorphism} $u\colon \phi \rightarrow \psi$ is an 
element of $L\{\tau\}$ such that $u \circ \phi_a = \psi_a \circ u$ for all  
$a \in A$. Hence an {\it endomorphism} of $\phi$ is an element of the
centralizer of $\phi(A)$ in $L\{\tau\}$, and $u$ is an {\it isomorphism}
if $u \in L^{\ast} \hookrightarrow L\{\tau\}$ is subject to $u \circ
\phi_a = \psi_a \circ u$.
\enddf

Define $\deg\colon  a \rightarrow \Bbb Z \cup \{-\infty\}$ and
$\deg_{\tau}\colon  L\{\tau\} \rightarrow \Bbb Z \cup \{-\infty\}$ by
$\deg(a) = \log_q|A/a|$ ($a \not= 0$; we write $A/a$ for $A/aA$),
$\deg (0) = -\infty$, and $\deg_{\tau}(f) =$ the well defined degree of $f$ as
a ``polynomial'' in $\tau$. It is an easy exercise in Dedekind rings to prove
the following

\th Proposition 1

 If $\phi$ is a Drinfeld module over $L$, there
exists a non-negative integer $r$  such that $\deg_{\tau}(\phi_a) 
= r \deg (a)$ for all $a \in A${\rm;} $r$ is called the {\it rank} ${\rm rk}(\phi)$
of $\phi$.
\endth

Obviously, ${\rm rk}(\phi) = 0$ means that $\phi = \gamma$, i.e., the 
$A$-module structure on ${\Bbb G_a}_{/L}$ is the tautological one.

\df Definition 4

Denote by ${\Cal M}^r(1)(L)$ the set of isomorphism
classes of Drinfeld modules of rank $r$ over $L$.
\enddf

\eg Example 3

 Let $A = \Bbb F_q[T]$ be as in Example 1 and let 
$K = \Bbb F_q(T)$ be its fraction field. 
Defining a Drinfeld module $\phi$ 
over $K$ or an extension field $L$ of $K$ is equivalent to specifying 
$\phi_T = T+g_1\tau + \cdots + g_r\tau^r \in L\{T\}$, where 
$g_r \not= 0$ and $r = {\rm rk}(\phi)$. In the special case where
$\phi_T = T+\tau$, $\phi$ is called the {\it Carlitz module}. Two such
Drinfeld modules $\phi$ and $\phi'$ are isomorphic over the algebraic
closure $L^{\alg}$ of $L$ if and only if  there is some $u \in {L^{\alg}}^{\,\ast}$
such that $g'_i = u^{q^i-1}g_i$ for all $i\ge 1$. Hence ${\Cal M}^r(1)
(L^{\alg})$ can be described (for $r \ge 1$) as an open dense 
subvariety of a weighted projective space of dimension $r-1$ over
$L^{\alg}$.
\endeg

\HH 4.2. Division points

\df Definition 5

 For $a \in A$ and a Drinfeld module $\phi$ over $L$,
write $_a\phi$ for the subscheme of {\it $a$-division points} of
${\Bbb G_a}_{/L}$ endowed with its structure of an $A$-module. Thus for any
$L$-algebra $M$,
 $$_a\phi(M) = \{x \in M\colon  \phi_a(x) = 0\}.$$
More generally, we put $_{\frak a}\phi = \bigcap\limits_{a\in {\frak a}} 
\phi_a$
for an arbitrary (not necessarily principal) ideal $\frak a$ of $A$. It is 
a finite flat group scheme of degree ${\rm rk}(\phi) \cdot \deg(\frak a)$,
whose structure is described in the next result.
\enddf

\th Proposition 2 {{\rm (\cite{Dr}, \cite{DH, I, Thm. 3.3 and Remark 3.4})}} 

 Let the Drinfeld module
$\phi$ over $L$ have rank $r \ge 1$.
\Roster 
\Item{(i)} If ${\rm char}_A(L) = \infty$, $_{\frak a}\phi$ is reduced for each
ideal $\frak a$ of $A$, and ${}_{\frak a}\phi(L^{\rm sep}) ={} _{\frak a}\phi
(L^{\alg})$ is isomorphic with $(A/\frak a)^r$ as an $A$-module. 
\Item{(ii)} If $\frak p = {\rm char}_A(L)$ is a maximal ideal, then there exists
an integer $h$, the {\it height} ${\rm ht}(\phi)$ of $\phi$, satisfying
$1 \le h \le r$, and such that $_{\frak a}\phi(L^{\alg})  \simeq 
(A/\frak a)^{r-h}$ whenever $\frak a$ is a power of $\frak p$, and 
$_{\frak a}\phi(L^{\alg})   \simeq (A/\frak a)^r$ if $(\frak a,\frak p)=1$.
\endRoster
\endth

The absolute Galois group $G_L$ of $L$ acts on $_{\frak a}\phi(L^{\rm sep})$
through $A$-linear automorphisms. Therefore, any Drinfeld module gives 
rise to Galois representations on its division points. These representations
tend to be ``as large as possible''.

The prototype of result is the following theorem, due to Carlitz
and Hayes \cite{H1}.

\th  Theorem 1

Let $A$ be the polynomial ring $\Bbb F_q[T]$ with field of
fractions $K$. Let $\rho\colon  A \rightarrow K\{\tau\}$ be
the Carlitz
module, $\rho_T = T + \tau$. For any non-constant monic polynomial
$a \in A$, let $K(a) := K(_a\rho(K^{\alg}))$ be the field extension
generated by the $a$-division points.
\Roster 
\Item{(i)} $K(a)/K$ is abelian with group $(A/a)^{\ast}$. If $\sigma_b$ is the
automorphism corresponding to the residue class of $b   \bmod a$ and
$x \in {}_a\rho(K^{\alg})$ then $\sigma_b(x) = \rho_b(x)$.
\Item{(ii)} If $(a) = \frak p^t$ is primary with some prime ideal $\frak p$ then
$K(a)/K$ is completely ramified at $\frak p$ and unramified at the
other finite primes.
\Item{(iii)} If $(a) = \prod a_i$ ($1 \le i \le s$) with primary and mutually
coprime $a_i$, the fields $K(a_i)$ are mutually linearly disjoint and 
$K = \otimes_{i \le i \le s} K(a_i)$.
\Item{(iv)} Let $K_+(a)$ be the fixed field of $\Bbb F_q^{\,\ast} \hookrightarrow
(A/a)^{\ast}$. Then $\infty$ is completely split in $K_+(a)/K$ and completely 
ramified in $K(a)/K_+(a)$.
\Item{(v)} Let $\frak p$ be a prime ideal generated by the monic polynomial
$\pi \in A$ and coprime with $a$. Under the identification
${\rm Gal}(K(a)/K) = (A/a)^{\ast}$,
the Frobenius element ${\rm Frob}_{\frak p}$ equals the residue class of
$\pi   \bmod a$.
\endRoster
\endth

Letting $a \rightarrow \infty$ with respect to divisibility, we obtain
the field $K(\infty)$ generated over $K$ by all the division points of
$\rho$, with group ${\rm Gal}(K(\infty)/K) = 
\inlim_{a}\,
(A/a)^{\ast}$, which almost agrees with the group of finite idele classes
of $K$. It turns out that $K(\infty)$ is the maximal abelian extension of
$K$ that is tamely ramified at $\infty$, i.e., we get a constructive version
of the class field theory of $K$. Hence the theorem may be seen both as a
global variant of Lubin--Tate's theory and as an analogue in characteristic
$p$ of the Kronecker--Weber theorem on cyclotomic extensions of $\Bbb Q$.

There are vast generalizations into two directions:
\Roster
\Item{(a)} abelian class field theory of arbitrary global function fields
$K = {\rm Frac}(A)$, where $A$ is a Drinfeld ring. 
\Item{(b)} systems of nonabelian Galois representations derived from 
Drinfeld modules.
\endRoster 

As to (a), the first problem is to find the proper analogue of the 
Carlitz module for an arbitrary Drinfeld ring $A$. As will result e.g. from
Theorem 2 (see also (4.3.4)), the isomorphism classes of rank-one Drinfeld
modules over the algebraic closure $K^{\alg}$ of $K$ correspond 
bijectively to the (finite!) class group ${\rm Pic}(A)$ of $A$. Moreover,
these Drinfeld modules $\rho^{(\frak a)}$ ($\frak a \in {\rm Pic}(A)$) may be
defined with coefficients in the ring ${\Cal O}_{H_+}$ of $A$-integers of
a certain abelian extension $H_+$ of $K$, and such that the leading
coefficients of all $\rho_a^{(\frak a)}$ are units of ${\Cal O}_{H_+}$.
Using these data along with the identification of $H_+$ in the 
dictionary of class field theory yields a generalization of Theorem 1 to the
case of arbitrary $A$. In particular, we again find an explicit construction
of the class fields of $K$ (subject to a tameness condition at $\infty$).
However, in view of class number problems, the theory (due to D. Hayes 
\cite{H2}, and superbly presented in \cite{Go2, Ch.VII}) has more of the
flavour of complex multiplication theory than of classical cyclotomic
theory.

Generalization (b) is as follows. Suppose that $L$ is a finite extension of
$K = {\rm Frac}(A)$, where $A$ is a general Drinfeld ring, and let the 
Drinfeld module $\phi$ over $L$ have rank $r$. For each power $\frak p^t$
of a prime $\frak p$ of $A$, $G_L = {\rm Gal}(L^{\rm sep}/L)$ acts on
$_{\frak p^t}\phi   \simeq (A/\frak p^t)^r$. We thus get an action of $G_L$ on the
$\frak p$-adic {\it Tate module} $T_{\frak p}(\phi)   \simeq (A_{\frak p})^r$ of
$\phi$ (see \cite{DH, I sect. 4}. Here of course $A_{\frak p} =
\prlim \,  A/\frak p^t$ is the $\frak p$-adic 
completion of $A$ with
field of fractions $K_{\frak p}$. Let on the other hand ${\rm End}(\phi)$ be
the endomorphism ring of $\phi$, which also acts on $T_{\frak p}(\phi)$. 
It is straightforward to show that 
(i) ${\rm End}(\phi)$ acts faithfully
and (ii) the two actions commute. In other words, we get an inclusion
 $$i\colon  {\rm End}(\phi)\otimes_{A} A_{\frak p}\hookrightarrow {\rm End}_{G_L}
  (T_{\frak p}(\phi))\leqno{(4.2.1)}$$
of finitely generated free $A_{\frak p}$-modules. The plain analogue of the
classical Tate conjecture for abelian varieties, proved 1983 by Faltings, 
suggests that $i$ is in fact bijective. This has been shown by Taguchi
\cite{Tag} and Tamagawa. Taking ${\rm End}(T_{\frak p}(\phi))   \simeq 
{\rm Mat}(r,A_{\frak p})$ and the known structure of subalgebras of
matrix algebras over a field into account, this means that the subalgebra
$$K_{\frak p}[G_L] \hookrightarrow {\rm End}(T_{\frak p}(\phi) 
\otimes_{A_{\frak p}} K_{\frak p})   \simeq {\rm Mat}(r,K_{\frak p})$$ 
generated
by the Galois operators is as large as possible. A much stronger statement
is obtained by R. Pink \cite{P1, Thm. 0.2}, who shows that the image of $G_L$
in ${\rm Aut}(T_{\frak p}(\phi))$ has finite index in the centralizer
group of ${\rm End}(\phi) \otimes A_{\frak p}$. Hence if e.g. $\phi$ has no 
``complex multiplications'' over $L^{\alg}$ (i.e., ${\rm End}_
{L^{\alg}}(\phi) = A$; this is the generic case for a Drinfeld module
in characteristic $\infty$), then the image of $G_L$ has finite index
in ${\rm Aut}(T_{\frak p}(\phi))   \simeq GL(r,A_{\frak p})$. This is
quite satisfactory, on the one hand, since we may use the Drinfeld module
$\phi$ to construct large nonabelian Galois extensions of $L$ with
prescribed ramification properties. On the other hand, the important
(and difficult) problem of estimating the index in question remains. 

\HH 4.3. Weierstrass theory

 Let $A$ be a Drinfeld ring with field of fractions  $K$, whose completion at $\infty$ is denoted by $K_{\infty}$.
We normalize the corresponding absolute value 
$|\,\,\, | = | \,\,\, |_{\infty}$ as $|a| = |A/a|$ for 
$0 \not= a\in A$ and let $C_{\infty}$ be the completed algebraic closure
of $K_{\infty}$, i.e., the completion of the algebraic closure
$K^{\alg}_{\infty}$ with respect to the unique extension of $|\,\,\,|$ to
$K^{\alg}_{\infty}$. By Krasner's theorem, $C_{\infty}$ is again
algebraically closed (\cite{BGS, p. 146}, where also other facts on
function theory in $C_{\infty}$ may be found). 
It is customary to indicate the strong
analogies between $A,K,K_{\infty},C_{\infty},\ldots$ and
$\Bbb Z,\Bbb Q,\Bbb R,\Bbb C,\ldots$, e.g. $A$ is a discrete and
cocompact subring of $K_{\infty}$. But note that $C_{\infty}$ fails
to be locally compact since $|C_{\infty}:K_{\infty}| = \infty$.

\df  Definition 6

 A {\it lattice of rank} $r$ (an $r$-lattice in brief)
in $C_{\infty}$ is a finitely generated (hence projective) discrete
$A$-submodule $\Lambda$ of $C_{\infty}$ of projective rank $r$, where the
discreteness means that $\Lambda$ has finite intersection with each
ball in $C_{\infty}$. The {\it lattice function} $e_{\Lambda}\colon 
C_{\infty} \rightarrow C_{\infty}$ of $\Lambda$ is defined as the
product
 $$e_{\Lambda}(z) = z \prod_{0 \not= \lambda \in \Lambda}
  (1-z/\lambda).\leqno{(4.3.1)}$$
It is entire (defined through an everywhere convergent power series),
$\Lambda$-periodic and $\Bbb F_q$-linear. For a non-zero $a\in A$ 
consider the diagram
 $$
\CD  
0 @>>> \Lambda  @>>>
           C_{\infty} @>{e_{\Lambda}}>> 
      C_{\infty} @>>> 0\\
   @. @V a VV @V a VV @V \phi_a^{\Lambda} VV \\
 0 @>>>  \Lambda 
    @>>>      C_{\infty} @>{e_{\Lambda}}>>   
      C_{\infty} @>>>  0 
\endCD \quad 
 \leqno{(4.3.2)}$$
with exact lines, where the left and middle arrows are multiplications by 
$a$ and $\phi_a^{\Lambda}$ is defined through commutativity. It is easy
to verify that 
\Roster
\Item{(i)} $\phi_a^{\Lambda} \in C_{\infty}\{\tau\}$, 

\Item{(ii)} $\deg_{\tau}(\phi_a^{\Lambda}) = r\cdot \deg(a)$, 

\Item{(iii)} $a \mapsto \phi_a^{\Lambda}$ is a ring homomorphism 
$\phi^{\Lambda}\colon  A \rightarrow C_{\infty}\{\tau\}$, in fact, a Drinfeld
module of rank $r$. Moreover, all the Drinfeld modules over $C_{\infty}$ are
so obtained.
\endRoster

\th  Theorem 2 {{\rm (Drinfeld \cite{Dr, Prop. 3.1})}}

\Roster
\Item{(i)}  Each rank-$r$ Drinfeld
module $\phi$ over $C_{\infty}$ comes via $\Lambda \mapsto \phi^{\Lambda}$
from some $r$-lattice $\Lambda$ in $C_{\infty}$.
\Item{(ii)}  Two Drinfeld modules $\phi^{\Lambda}$, $\phi^{\Lambda'}$ are 
isomorphic if and only if  there exists $0 \not= c \in C_{\infty}$ such that
$\Lambda' = c\cdot \Lambda$.
\endRoster 
\endth 

We may thus describe ${\Cal M}^r(1)(C_{\infty})$ (see Definition 4) as the 
space of $r$-lattices modulo similarities, i.e., as some generalized upper
half-plane modulo the action of an arithmetic group. Let us make this more
precise. 

\df Definition 7

 For $r \ge 1$ let $\Bbb P^{r-1}(C_{\infty})$ be
the $C_{\infty}$-points of projective $r-1$-space and $\Omega^r :=
\Bbb P^{r-1}(C_{\infty})-\bigcup H(C_{\infty})$, where $H$ runs through the
$K_{\infty}$-rational hyperplanes of $\Bbb P^{r-1}$. That is,
$ \Underline{\omega} = (\omega_1: \ldots :\omega_r)$ belongs to
{\it Drinfeld's half-plane} $\Omega^r$ if and only if  there is no non-trivial
relation $\sum a_i\omega_i = 0$ with coefficients $a_i \in K_{\infty}$.
\enddf

Both point sets $\Bbb P^{r-1}(C_{\infty})$ and $\Omega^r$ carry structures
 of analytic spaces over $C_{\infty}$ (even over $K_{\infty}$), and so we
can speak of holomorphic functions on $\Omega^r$. We will not give the
details (see for example \cite{GPRV, in particular lecture 6}); suffice it 
to say that locally uniform limits of rational functions (e.g. Eisenstein
series, see below) will be holomorphic.

Suppose for the moment that the class number $h(A) = |{\rm Pic}(A)|$ of
$A$ equals one, i.e., $A$ is a principal ideal domain. Then each
$r$-lattice $\Lambda$ in $C_{\infty}$ is free, 
$\Lambda = \sum_{1 \le i \le r} A\omega_i$, and the discreteness of
$\Lambda$ is equivalent with $ \Underline{\omega} := 
(\omega_1: \ldots :\omega_r)$ belonging to $\Omega^r \hookrightarrow 
\Bbb P^{r-1}(C_{\infty})$. Further, two points $ \Underline{\omega}$ and
$ \Underline{\omega}'$ describe similar lattices (and therefore isomorphic
Drinfeld modules) if and only if  they are conjugate under $\Gamma := GL(r,A)$,
which acts on $\Bbb P^{r-1}(C_{\infty})$ and its subspace $\Omega^r$.
Therefore, we get a canonical bijection
 $$\Gamma\setminus \Omega^r  \iss 
  {\Cal M}^r(1)(C_{\infty})\leqno{(4.3.3)}$$
from the quotient space $\Gamma \setminus \Omega^r$ to the set of
isomorphism classes ${\Cal M}^r(1)(C_{\infty})$.

In the general case of arbitrary $h(A) \in \Bbb N$, we let 
$\Gamma_i := GL(Y_i) \hookrightarrow GL (r,k)$, where
$Y_i \hookrightarrow K^r$ ($1 \le i \le h(A)$) runs through 
representatives of the $h(A)$ isomorphism classes of projective 
$A$-modules of rank $r$. In a similar fashion (see e.g. \cite{G1, II sect.1},
\cite{G3}), we get a bijection
 $${\bigcupp}_{1 \le i \le h(A)} 
 \Gamma_i \setminus 
  \Omega^r  \iss {\Cal M}^r(1)(C_{\infty}),
 \leqno{(4.3.4)}$$
 which can be made independent of choices if we use the canonical
adelic description of the $Y_i$.

\eg Example 4

 If $r=2$ then $\Omega = \Omega^2 = \Bbb P^1(C_{\infty})
- - \Bbb P^1(K_{\infty}) = C_{\infty}-K_{\infty}$, which rather corresponds
to $\Bbb C-\Bbb R = H^+ \bigcup H^-$ (upper and lower complex half-planes)
than to $H^+$ alone. The group $\Gamma := GL(2,A)$ acts
via ${a\,b\choose c\,d} (z) = \frac{az+b}{cz+d}$, and thus gives rise
to {\it Drinfeld modular forms} on $\Omega$ (see \cite{G1}). Suppose moreover
that $A = \Bbb F_q[T]$ as in Examples 1 and 3. We define {\it ad hoc}
a modular form of weight $k$ for $\Gamma$ as a holomorphic function
$f\colon  \Omega \rightarrow C_{\infty}$ that satisfies 
\Roster
\Item{(i)} $f\left(\frac{az+b}{cz+d}\right) = (cz+d)^k f(z)$ for ${a\,b\choose
c\,d}\in \Gamma$ and 
\endRoster
\smallskip
\Roster
\Item{(ii)} $f(z)$ is bounded on the subspace $\{z \in \Omega :  
{\rm inf}_{x \in K_{\infty}} |z-x| > 1\}$ of $\Omega$.
\endRoster 
Further, we put $M_k$ for the $C_{\infty}$-vector space of modular forms
of weight $k$. (In the special case under consideration, (ii) is
equivalent to the usual ``holomorphy at cusps'' condition. For more
general groups $\Gamma$, e.g. congruence subgroups of
$GL(2,A)$, general Drinfeld rings $A$, and higher ranks $r \ge 2$,
condition (ii) is considerably more costly to state, see \cite{G1}.) Let
 $$E_k(z) := \sum_{(0,0)\not=(a,b) \in A \times A} 
  \frac{1}{(az+b)^k} \leqno{(4.3.5)}$$
be the {\it Eisenstein series} of weight $k$. Due to the non-archimedean
situation, the sum converges for $k \ge 1$ and yields a modular form
$0 \not= E_k \in M_k$ if $k \equiv 0~(q-1)$. Moreover, the various $M_k$
are linearly independent and 
 $$M(\Gamma) := \bigoplus_{k \ge 0} M_k = C_{\infty}[E_{q-1},E_{q^2-1}]
  \leqno{(4.3.6)}$$
is a polynomial ring in the two algebraically independent Eisenstein
series of weights $q-1$ and $q^2-1$. There is an {\it a priori} different
method of constructing modular forms via Drinfeld modules. With each
$z \in \Omega$, associate the 2-lattice $\Lambda_z := Az+A \hookrightarrow
C_{\infty}$ and the Drinfeld module $\phi^{(z)} = \phi^{(\Lambda_z)}$.
Writing $\phi_T^{(z)} = T+g(z)\tau+\Delta(z)\tau^2$, the 
coefficients $g$ and $\Delta$ become functions in $z$, in fact, modular
forms of respective weights $q-1$ and $q^2-1$. We have (\cite{Go1}, 
\cite{G1, II 2.10})
 $$g = (T^g-T)E_{q-1},\: \Delta=(T^{q^2}-T)E_{q^2-1}+(T^{q^2}-T^q)
  E_{q-1}^{q+1}.\leqno{(4.3.7)}$$
The crucial fact is that $\Delta(z) \not= 0$ for $z \in \Omega$, but 
$\Delta$ vanishes ``at infinity''. Letting $j(z) := g(z)^{q+1}/\Delta(z)$
(which is a function on $\Omega$ invariant under $\Gamma$), the 
considerations of Example 3 show that $j$ is a complete invariant for
Drinfeld modules of rank two. Therefore, the composite map
 $$j\colon  \Gamma \setminus \Omega \iss 
  {\Cal M}^2(1)(C_{\infty})  \iss C_{\infty}
  \leqno{(4.3.8)}$$
is bijective, in fact, biholomorphic. 

\HH 4.4. Moduli schemes

 We want to give a similar description of
${\Cal M}^r(1)(C_{\infty})$ for $r \ge 2$ and arbitrary $A$, that is,
to convert (4.3.4) into an isomorphism of analytic spaces. One proceeds as
follows (see \cite{Dr}, \cite{DH}, \cite{G3}): 

\noindent{(a)} Generalize the notion of ``Drinfeld $A$-module over an $A$-field $L$''
to ``Drinfeld $A$-module over an $A$-scheme $S \rightarrow {\rm Spec}~A$''.
This is quite straightforward. Intuitively, a Drinfeld module over $S$
is a continuously varying family of Drinfeld modules over the residue
fields of $S$.

\noindent{(b)} Consider the functor on $A$-schemes:
 $${\Cal M}^r\colon  S \longmapsto \left\{
\aligned 
 &\text{\rm isomorphism classes of rank-$r$}\\
  &\text{\rm Drinfeld modules over $S$}
\endaligned \right\}.$$
The naive initial question is to represent this functor by an
$S$-scheme $M^r(1)$. This is impossible in view of the existence of
automorphisms of Drinfeld modules even over algebraically closed
$A$-fields.

\noindent{(c)} As a remedy, introduce rigidifying level structures on 
Drinfeld modules.
Fix some ideal $0 \not= \frak n$ of $A$. An {\it $\frak n$-level structure}
on the Drinfeld module $\phi$ over the $A$-field $L$ whose $A$-characteristic
doesn't divide $\frak n$ is the choice of an isomorphism of $A$-modules
 $$\alpha\colon  (A/\frak n)^r \iss 
 _{\frak n}\phi(L)$$
(compare Proposition 2). Appropriate modifications apply to the cases where
${\rm char}_A(L)$ divides $\frak n$ and where the definition field $L$ is
replaced by an $A$-scheme $S$. Let ${\Cal M}^r(\frak n)$ be the functor
$${\Cal M}^r(\frak n)\colon  S \longmapsto \left\{
\aligned
&\text{\rm isomorphism classes of rank-$r$}\\
 &\text{\rm Drinfeld modules over $S$ endowed} \\
  &\text{\rm with an $\frak n$-level structure}
\endaligned\right\}.$$

\th   Theorem 3 {{\rm (Drinfeld \cite{Dr, Cor. to Prop. 5.4})}}

  Suppose that $\frak n$ is divisible
by at least two different prime ideals. Then ${\Cal M}^r(\frak n)$ is
representable by a smooth affine $A$-scheme $M^r(\frak n)$ of
relative dimension $r-1$.
\endth

In other words, the scheme $M^r(\frak n)$ carries a ``tautological'' Drinfeld
module $\phi$ of rank $r$ endowed with a level-$\frak n$ structure such
that pull-back induces for each $A$-scheme $S$ a bijection
$$
M^r(\frak n)(S) = \{\text{\rm morphisms }  (S,M^r(\frak n)) \text{\rm\} }  
 {}   \iss  {\Cal M}^r(\frak n)(S), \quad 
    f {}\longmapsto f^{\ast}(\phi).   
\leqno{(4.4.1)}
 $$ 
$M^r(\frak n)$ is called the (fine) {\it moduli scheme} for the moduli
problem ${\Cal M}^r(\frak n)$. Now the finite group $G(\frak n) :=
GL(r,A/\frak n)$ acts on ${\Cal M}^r(\frak n)$ by permutations of
the level structures. By functoriality, it also acts on $M^r(\frak n)$.
We let $M^r(1)$ be the quotient of $M^r(\frak n)$ by $G(\frak n)$ (which
does not depend on the choice of $\frak n$). It has the property that at 
least its $L$-valued points for algebraically closed $A$-fields $L$ 
correspond bijectively and functorially to ${\Cal M}^r(1)(L)$. It is
therefore called a {\it coarse moduli scheme} for ${\Cal M}^r(1)$.
Combining the above with (4.3.4) yields a bijection
 $$\bigcupp_{1 \le i \le h(A)} \Gamma_i\setminus 
   \Omega^r \iss  M^r(1)(C_{\infty}),\leqno{(4.4.2)}$$ 
which even is an isomorphism of the underlying analytic 
spaces \cite{Dr, Prop. 6.6}. The most simple special case is the one dealt with
in Example 4, where $M^2(1) = {\Bbb A^1}_{/A}$, the affine line over $A$. 

\HH 4.5. Compactification

 It is a fundamental question to construct
and study a ``compactification'' of the affine $A$-scheme $M^r(\frak n)$,
relevant for example for the Langlands conjectures over $K$, the cohomology of
arithmetic subgroups of $GL(r,A)$, or the $K$-theory of $A$ and $K$.
This means that we are seeking  a proper $A$-scheme $ \Overline{M}^r(\frak n)$
with an $A$-embedding $M^r(\frak n) \hookrightarrow  \Overline{M}^r(\frak n)$
as an open dense subscheme, and which behaves functorially with respect
to the forgetful morphisms $M^r(\frak n) \rightarrow M^r(\frak m)$
if $\frak m$ is a divisor of $\frak n$. For many purposes it suffices
to solve the apparently easier problem of constructing similar
compactifications of the generic fiber $M^r(\frak n) \times_A K$
or even of $M^r(\frak n) \times_A C_{\infty}$. Note that varieties
over $C_{\infty}$ may be studied by analytic means, using the GAGA 
principle.

There are presently three approaches towards the problem of compactification:

\noindent(a) a (sketchy) construction of the present author \cite{G2} of a 
compactification $ \Overline{M}_{\Gamma}$ of $M_{\Gamma}$, the
$C_{\infty}$-va\-ri\-e\-ty corresponding to an arithmetic subgroup $\Gamma$
of $GL(r,A)$ (see (4.3.4) and (4.4.2)). We will return to this
below;

\noindent(b) an analytic compactification similar to (a), restricted to the case
of a polynomial ring $A = \Bbb F_q[T]$, but with the advantage of presenting
complete proofs, by M. M. Kapranov \cite{K};

\noindent(c) R. Pink's idea of a modular compactification of $M^r(\frak n)$ over $A$
through a generalization of the underlying moduli problem \cite{P2}.

Approaches (a) and (b) agree essentially in their common domain, up to notation
and some other choices. Let us briefly describe how one proceeds in (a).
Since there is nothing  to show for $r=1$, we suppose that $r \ge 2$.

We let $A$ be any Drinfeld ring. If $\Gamma$ is a subgroup of $GL(r,K)$ 
commensurable with $GL(r,A)$ (we call such $\Gamma$ 
{\it arithmetic subgroups}), the point set $\Gamma \setminus \Omega$ is the
set of $C_{\infty}$-points of an affine variety $M_{\Gamma}$ over $C_{\infty}$,
as results from the discussion of subsection 4.4. If $\Gamma$ is the congruence
subgroup $\Gamma(\frak n) = \{\gamma \in GL(r,A)\colon  \gamma \equiv 1 \bmod 
\frak n\}$, then $M_{\Gamma}$ is one of the irreducible components of 
$M^r(\frak n) \times_A C_{\infty}$.

\df  Definition 8

 For $ \Underline{\omega} = (\omega_1,:\ldots :\omega_r)
\in \Bbb P^{r-1}(C_{\infty})$ put 
$$
r( \Underline{\omega}) := \dim_K
(K\omega_1 + \cdots + K\omega_r)\quad{\text{and}}\quad r_{\infty}( \Underline{\omega})
:= \dim_{K_{\infty}} (K_{\infty}\omega_1+ \cdots + K_{\infty} \omega_r).$$
Then $1 \le r_{\infty}( \Underline{\omega}) \le r( \Underline{\omega})
\le r$ and $\Omega^r = \{ \Underline{\omega}~|~r_{\infty}( \Underline
{\omega})=r\}$. More generally, for $1 \le i \le r$ let
 $$\Omega^{r,i} := \{ \Underline{\omega}\colon  r_{\infty}( \Underline{\omega})
  = r( \Underline{\omega})= i\}.$$
Then $\Omega^{r,i}=\bigcupp \Omega_V$, where $V$ runs 
through the $K$-subspaces of dimension $i$ of $K^r$ and $\Omega_V$ is 
constructed from $V$ in a similar way as is $\Omega^r = \Omega_{K^r}$
from $C_{\infty}^r = (K^r) \otimes C_{\infty}$. That is,
$\Omega_V = \{ \Underline{\omega} \in \Bbb P(V \otimes C_{\infty})
\hookrightarrow \Bbb P^{r-1}(C_{\infty})\colon  r_{\infty}( \Underline{\omega})
= r( \Underline{\omega}) =i\}$, which has a natural structure as analytic
space of dimension $\dim(V)-1$ isomorphic with $\Omega^{\dim(V)}$.
Finally, we let $ \Overline{\Omega}^r := \{ \Underline{\omega}\colon 
r_{\infty}( \Underline{\omega}) = r( \Underline{\omega})\} =
{\bigcupp}_{1 \le i \le r} \Omega^{r,i}$. 
\enddf 

$ \Overline{\Omega}^r$ along with its stratification through the
$\Omega^{r,i}$ is stable under $GL(r,K)$, so this also holds for
the arithmetic group $\Gamma$ in question. The quotient $\Gamma
\setminus  \Overline{\Omega}^r$ turns out to be the $C_{\infty}$-points of
the wanted compactification $ \Overline{M}_{\Gamma}$.

\df Definition 9

 Let $P_i \hookrightarrow G := GL(r)$ be the
maximal parabolic subgroup of matrices with first $i$ columns being zero. 
Let $H_i$ be the obvious factor group isomorphic  $GL(r-i)$.  
Then $P_i(K)$ acts via $H_i(K)$ on $K^{r-i}$ and thus on $\Omega^{r-i}$.
From $$G(K)/P_i(K)  \iss 
\{\text{\rm subspaces $V$ of dimension $r-i$ of $K^r$}\}$$ we get bijections
 $$
 \aligned
&G(K) \times_{P_i(K)} \Omega^{r-i}   \iss 
  \Omega^{r,r-i}, \\
 & (g,\omega_{i+1}: \ldots :\omega_r)  \longmapsto  
  (0: \cdots :0:\omega_{i+1}: \ldots :\omega_r)g^{-1}
 \endaligned
 \leqno{(4.5.1)}$$
and
 $$\Gamma\setminus \Omega^{r,r-i} \iss
{\bigcupp}_{g\in \Gamma\setminus G(K)/P_i(K)}
  \Gamma(i,g)\setminus \Omega^{r-i},\leqno{(4.5.2)}$$
where $\Gamma(i,g) := P_i \cap g^{-1}\Gamma_g$, and the double
quotient $\Gamma\setminus G(K)/P_i(K)$ is finite by elementary lattice
theory. Note that the image of $\Gamma(i,g)$ in $H_i(K)$ (the
group that effectively acts on $\Omega^{r-i}$) is again an arithmetic 
subgroup of $H_i(K) = GL(r-i,K)$, and so the right hand side
of (4.5.2) is the disjoint union of analytic spaces of the same type
$\Gamma'\setminus \Omega^{r'}$.
\enddf

\eg  Example 5

 Let $\Gamma = \Gamma(1) = GL(r,A)$ and $i=1$. Then
$\Gamma\setminus G(K)/P_1(K)$ equals the set of isomorphism classes
of projective $A$-modules of rank $r-1$, which in turn (through the
determinant map) is in one-to-one correspondence with the class group
${\rm Pic}(A)$.
\endeg

Let $F_V$ be the image of $\Omega_V$ in $\Gamma\setminus \Overline{\Omega}^r$.
The different analytic spaces $F_V$, corresponding to locally closed 
subvarieties of $ \Overline{M}_{\Gamma}$, are glued together in such a way
that $F_U$ lies in the Zariski closure $ \Overline{F}_V$ of $F_V$ if and only if  $U$
is $\Gamma$-conjugate to a $K$-subspace of $V$. Taking into account
that $F_V   \simeq \Gamma'\setminus \Omega^{\dim(V)} = M_{\Gamma'}(C_{\infty})$
for some arithmetic subgroup $\Gamma'$ of $GL(\dim(V),K)$,
$ \Overline{F}_V$ corresponds to the compactification $ \Overline{M}_{\Gamma'}$
of $M_{\Gamma'}$.

The details of the gluing procedure are quite technical and complicated
and cannot be presented here (see \cite{G2} and \cite{K} for some special
cases). Suffice it to say that for each boundary component $F_V$ of
codimension one, a vertical coordinate $t_V$ may be specified such that 
$F_V$ is locally given by $t_V=0$. The result (we refrain from stating a
``theorem'' since proofs of the assertions below in full strength and
generality are published neither in \cite{G2} nor in \cite{K}) will be
a normal projective $C_{\infty}$-variety $ \Overline{M}_{\Gamma}$ provided
with an open dense embedding $i\colon  M_{\Gamma} \hookrightarrow 
 \Overline{M}_{\Gamma}$ with the following properties:
 
 \bitem $ \Overline{M}_{\Gamma}(C_{\infty}) = \Gamma\setminus
 \Overline{\Omega}^r$, and the inclusion $\Gamma \setminus \Omega^r
\hookrightarrow \Gamma\setminus  \Overline{\Omega}^r$ corresponds to $i$;
 \bitem $ \Overline{M}_{\Gamma}$ is defined over the same finite abelian
extension of $K$ as is $M_{\Gamma}$;
 \bitem for $\Gamma' \hookrightarrow \Gamma$, the natural map 
$M_{\Gamma'} \rightarrow M_{\Gamma}$ extends to 
$ \Overline{M}_{\Gamma'} \rightarrow  \Overline{M}_{\Gamma}$;
 \bitem the $F_V$ correspond to locally closed subvarieties, and 
$ \Overline{F}_V = \cup F_U$, where $U$ runs through the $K$-subspaces of
$V$ contained up to the action of $\Gamma$ in $V$;
 \bitem $ \Overline{M}_{\Gamma}$ is ``virtually non-singular'', i.e., $\Gamma$
contains a subgroup $\Gamma'$ of finite index such that 
$ \Overline{M}_{\Gamma'}$ is non-singular; in that case, the boundary components
of codimension one present normal crossings.

\smallskip \par 
 
Now suppose that $ \Overline{M}_{\Gamma}$ is non-singular and that 
$x \in  \Overline{M}_{\Gamma}(C_{\infty}) = \bigcup_{1\le i \le r}
\Omega^{r,i}$ belongs to $\Omega^{r,1}$. Then we can find a sequence
$\{x\} = X_0 \subset \cdots X_i \cdots \subset X_{r-1} = 
 \Overline{M}_{\Gamma}$ of smooth subvarieties $X_i =  \Overline{F}_{V_i}$
of dimension $i$. Any holomorphic function around $x$ (or more generally,
any modular form for $\Gamma$) may thus be expanded as a series in
$t_V$ with coefficients in the function field of $ \Overline{F}_{V_{r-1}}$,
etc. Hence $ \Overline{M}_{\Gamma}$ (or rather its completion at the $X_i$) 
may be described through $(r-1)$-dimensional local fields with residue
field $C_{\infty}$. The expansion of some standard modular forms can be
explicitly calculated, see \cite{G1, VI} for the case of $r=2$. In the last
section we shall present at least the vanishing orders of some of these
forms.

\eg   Example 6

 Let $A$ be the polynomial ring $\Bbb F_q[T]$ and 
$\Gamma = GL(r,A)$. As results from Example 3, (4.3.3) and (4.4.2),
 $$M_{\Gamma}(C_{\infty}) = M^r(1)(C_{\infty}) = \{(g_1,\ldots,g_r) \in
C_{\infty}^r\colon  g_r \not= 0\}/C_{\infty}^{\ast},$$
where $C_{\infty}^{\ast}$ acts diagonally through $c(g_1,\ldots,g_r) = 
(\ldots,c^{q^i-1}g_i,\ldots)$, which is the open subspace of weighted
projective space $\Bbb P^{r-1}(q-1,\ldots,q^r-1)$ with non-vanishing last
coordinate. The construction yields
 $$   \Overline{M}_{\Gamma}(C_{\infty})  =  \Bbb P^{r-1}(q-1,\ldots,q^r-1)
  (C_{\infty})
    =   {\bigcupp}_{1 \le i \le r} M^i(1)(C_{\infty}).$$
Its singularities are rather mild and may be removed upon replacing $\Gamma$
by a congruence subgroup.
\endeg 

\HH 4.6. Vanishing orders of modular forms

 In this final section 
we state some results about the vanishing orders of certain modular forms
along the boundary divisors of $ \Overline{M}_{\Gamma}$, in  the case where
$\Gamma$ is either $\Gamma(1) = GL(r,A)$ or a full congruence 
subgroup $\Gamma(\frak n)$ of $\Gamma(1)$. These are relevant for the
determination of $K$- and Chow groups, and for standard conjectures about
the arithmetic interpretation of partial zeta values.

In what follows, we suppose that $r \ge 2$, and put $z_i := 
\frac{\omega_i}{\omega_r}$ ($1 \le i \le r$) for the coordinates
$(\omega_1: \ldots :\omega_r)$ of $ \Underline{\omega} \in \Omega^r$.
Quite generally, $ \Underline{a} = (a_1,\ldots,a_r)$ denotes a vector with
$r$ components.

\df  Definition 10

 The {\it Eisenstein series} $E_k$ of weight $k$ on
$\Omega^r$ is defined as
 $$E_k( \Underline{\omega}) := \sum_{ \Underline{0} \not=  \Underline{a} \in A^r}
  \frac{1}{(a_1z_1+ \cdots + a_rz_r)^k}.$$
Similarly, we define for $ \Underline{u} \in \frak n^{-1} \times \cdots 
\times \frak n^{-1} \subset K^r$
 $$E_{k, \Underline{u}}( \Underline{\omega}) = \sum_{ \Underline{0} \not= 
   \Underline{a}\in K^r \atop  \Underline{a} \equiv  \Underline{u}   \bmod A^r}
  \frac{1}{(a_1z_1+\cdots + a_rz_r)^k}.$$
These are modular forms for $\Gamma(1)$ and $\Gamma(\frak n)$, respectively, 
that is, they are holomorphic, satisfy the obvious transformation values
under $\Gamma(1)$ (resp. $\Gamma(\frak n)$), and extend to sections of a 
line bundle on $ \Overline{M}_{\Gamma}$. As in Example 4, there is a second
type of modular forms coming directly from Drinfeld modules.
\enddf

\df  Definition 11

 For $ \Underline{\omega} \in \Omega^r$ write
$\Lambda_{ \Underline{\omega}} = Az_1+ \cdots + Az_r$ and 
$e_{ \Underline{\omega}}$, $\phi^{ \Underline{\omega}}$ for the lattice
function and Drinfeld module associated with 
$\Lambda_{ \Underline{\omega}}$, respectively.
If $a \in A$ has degree $d = \deg(a)$,
 $$\phi_a^{ \Underline{\omega}} = a+\sum_{1 \le i \le r\cdot d} 
  \ell_i(a, \Underline{\omega})\tau^i.$$
\endeg

The $\ell_i(a, \Underline{\omega})$ are modular forms of weight $q^i-1$ 
for $\Gamma$. This holds in particular for
 $$\Delta_a( \Underline{\omega}) := \ell_{rd}(a, \Underline{\omega}),$$
which has weight $q^{rd}-1$ and vanishes nowhere on $\Omega^r$. The functions
$g$ and $\Delta$ in Example 4 merely constitute a very special instance
of this construction. We further let, for $ \Underline{u} \in 
(\frak n^{-1})^r$,
 $$e_{ \Underline{u}}( \Underline{\omega}) := e_{ \Underline{\omega}}
  (u_1z_1+ \cdots + u_rz_r),$$
the $\frak n$-division point of type $ \Underline{u}$ of 
$\phi^{ \Underline{\omega}}$. If $ \Underline{u} \not\in A^r$, 
$e_{ \Underline{u}}( \Underline{\omega})$ vanishes nowhere on $\Omega^r$,
and it can be shown that in this case,
 $$e_{ \Underline{u}}^{-1} = E_{1, \Underline{u}}.\leqno{(4.6.1)}$$
We are interested in the behavior around the boundary of 
$ \Overline{M}_{\Gamma}$ of these forms. Let us first describe the set
$\{ \Overline{F}_V\}$ of boundary divisors, i.e., of irreducible components,
all of codimension one, of $ \Overline{M}_{\Gamma}-M_{\Gamma}$. For
$\Gamma = \Gamma(1) = GL(r,A)$, there is a natural bijection
 $$\{ \Overline{F}_V\}  \iss {\rm Pic}(A)
  \leqno{(4.6.2)}$$
described in detail in \cite{G1, VI 5.1}. It is induced from 
$V \mapsto \text{\rm inverse of } \Lambda^{r-1}(V\cap A^r)$. (Recall that
$V$ is a $K$-subspace of dimension $r-1$ of $K^r$, thus $V \cap A^r$ a
projective module of rank $r-1$, whose $(r-1)$-th exterior power 
$\Lambda^{r-1}(V\cap A^r)$ determines an element of ${\rm Pic}(A)$.)
We denote the component corresponding to the class $(\frak a)$ of an ideal
$\frak a$ by $ \Overline{F}_{(\frak a)}$. Similarly, the boundary divisors
of $ \Overline{M}_{\Gamma}$ for $\Gamma = \Gamma(\frak n)$ could be 
described via generalized class groups. We simply use (4.5.1) and
(4.5.2), which now give 
 $$\{ \Overline{F}_V\}  \iss \Gamma(\frak n)
  \setminus GL(r,K)/P_1(K).\leqno{(4.6.3)}$$
We denote the class of $\nu \in GL(r,K)$ by $[\nu]$. For the 
description of the behavior of our modular forms along the $ \Overline{F}_V$,
we need the partial zeta functions of $A$ and $K$. For more about these,
see \cite{W} and \cite{G1, III}.

\df  Definition 12

 We let
 $$\zeta_K(s) = \sum |\frak a|^{-s} = \frac{P(q^{-s})}{(1-q^{-s})(1-q^{1-s})}$$
be the zeta function of $K$ with numerator polynomial $P(X) \in \Bbb Z[X]$.
Here the sum is taken over the positive divisors $\frak a$ of $K$ (i.e., of the 
curve ${\Cal C}$ with function field $K$). Extending the sum only over 
divisors with support in ${\rm Spec}(A)$, we get
 $$\zeta_A(s) = \sum_{0 \not= \frak a \subset A~{\rm ideal}}
  |\frak a|^{-s} = \zeta_K(s)(1-q^{-d_{\infty}s}),$$
where $d_{\infty} = \deg_{\Bbb F_q}(\infty)$. For a class $\frak c \in
{\rm Pic}(A)$ we put 
 $$\zeta_{\frak c}(s) = \sum_{\frak a \in \frak c}|\frak a|^{-s}.$$
If finally $\frak n \subset K$ is a fractional $A$-ideal and $t \in K$, we
define 
 $$\zeta_{t  \bmod \frak n}(s) = \sum_{a\in K\atop a \equiv t  \bmod \frak n}
  |a|^{-s}.$$
Among the obvious distribution relations \cite{G1, III sect.1} between
these, we only mention
 $$\zeta_{(\frak n^{-1})}(s) = \frac{|\frak n|^s}{q-1} \zeta_{0  \bmod \frak n}
  (s).\leqno{(4.6.4)}$$
We are now in a position to state the following theorems, which may be
proved following the method of \cite{G1, VI}.

\th Theorem 4

  Let $a\in A$ be non-constant and $\frak c$ a class
in ${\rm Pic}(A)$. The modular form $\Delta_a$ for $GL(r,A)$ has
vanishing order 
 $${\rm ord}_{\frak c}(\Delta_a) = -(|a|^r-1)\zeta_{\frak c}(1-r)$$
at the boundary component $ \Overline{F}_{\frak c}$ corresponding to
$\frak c$.
\endth

\th Theorem 5
 
 Fix an ideal $\frak n$ of $A$ and $ \Underline{u}
\in K^r-A^r$ such that $ \Underline{u}\cdot \frak n \subset A^r$, 
and let $e_{ \Underline{u}}^{-1} = E_{1, \Underline{u}}$ be the modular
form for $\Gamma(\frak n)$ determined by these data. The vanishing order
${\rm ord}_{[\nu]}$ of $E_{1, \Underline{u}}( \Underline{\omega})$ at the
component corresponding to $\nu \in GL(r,K)$ {{\rm(}}see $(4.6.2)${{\rm)}} is
given as follows{{\rm:}} \ let $\pi_1\colon  K^r \rightarrow K$ be the projection
to the first coordinate and let $\frak a$ be the fractional ideal 
$\pi_1(A^r\cdot \nu)$. Write further $ \Underline{u} \cdot \nu =
(v_1,\ldots,v_r)$. Then  
 $${\rm ord}_{[\nu]} E_{1, \Underline{u}}(\omega) = \frac{|\frak n|^{r-1}}
  {|\frak a|^{r-1}} (\zeta_{v_1 \bmod \frak a}(1-r)-\zeta_{0 \bmod \frak a}
  (1-r)).$$
\endth

Note that the two theorems do not depend on the full strength of 
properties of $ \Overline{M}_{\Gamma}$ as stated without proofs in the
last section, but only on the {\it normality} of $ \Overline{M}_{\Gamma}$,
which is {\it proved} in $[K]$ for $A = \Bbb F_q[T]$, and whose 
generalization to arbitrary Drinfeld rings is straightforward (even
though technical).

 \Bib References

 \rf{BGS} S. Bosch and U. G\"untzer and R. Remmert, Non-archimedean
Analysis, Grundl. Math. Wiss. {  261}, Springer 1984.
 
\rf{DH} P. Deligne and D. Husem\"oller, Survey of Drinfeld modules,
Contemp. Math. {  67} (1987), 25--91.
 \rf{Dr} V.G. Drinfeld, Elliptic modules, Math. USSR-Sbornik
{  23} (1976), 561--592.

 \rf{G1} E.-U. Gekeler, Drinfeld modular curves, Lect. Notes
Math. {  1231}, Springer 1986.

 \rf{G2} E.-U. Gekeler, Satake compactification of Drinfeld
modular schemes, Proc. Conf. $P$-adic Analysis, Hengelhoef 1986.

 \rf{G3} E.-U. Gekeler, Moduli for Drinfeld modules, in:
``The Arithmetic of Function Fields'', D. Goss et al. (eds.), Walter
de Gruyter 1992.

 \rf{GPRV} Drinfeld Modules, Modular Schemes and Applications,
E.-U. Gekeler et al. (eds.), World Scientific 1997.

 \rf{Go1} D. Goss, Modular forms for $\Bbb F_r[T]$, J. reine
angew. Math. {  317} (1980), 3--38.
 
\rf{Go2} D. Goss, Basic Structures of Function Field Arithmetic,
Ergeb. Math. {  35}, Springer 1996.
 
\rf{H1} D. Hayes, Explicit class field theory for rational
function fields, Trans. AMS {  189} (1974), 77--91.

 \rf{H2} D. Hayes, Explicit class field theory in global 
function fields, in: G.C. Rota (ed.), Studies in Algebra and Number
Theory, Academic Press 1979.
 
\rf{K} M.M. Kapranov, On cuspidal divisors on the modular
varieties of elliptic modules, Math. USSR Izv. {  30} (1988), 533--547.
 
\rf{P1} R. Pink, The Mumford-Tate conjecture for Drinfeld
modules, Publ. RIMS {  33} (1997), 393--425.

 \rf{P2} R. Pink, unpublished notes.
 
\rf{Tag} Y. Taguchi, The Tate conjecture for $t$-motives,
Proc. AMS {  123} (1995), 3285--3287.

 \rf{Tam} A. Tamagawa, The Tate conjecture for $A$-premotives,
Preprint 1994.

 \rf{Tat} J. Tate, Les Conjectures de Stark sur les Functions
$L$ d'Artin en $s=0$, Prog. Math. {  47}, Birkh\"auser 1984.
 
\rf{W} A. Weil, Basic Number Theory, Springer 1967.
 \endBib

\Coordinates 

FR 6.1 - Mathematik \ 
Universit\"at des Saarlandes

Postfach 15 11 50 \ 
D-66041 Saarbr\"ucken \ Germany 

E-mail: gekeler\@math.uni-sb.de
\endCoordinates

\vfill
\pagebreak

\end

%% file: m3-macs.tex
%%%% prevent double loading:
\expandafter\ifx\csname mthreemacsloaded\endcsname\relax\else \fi

\magnification1100
\input amstex

%%% Hack of Plain TeX correction and style macros 
%%% written by Walter Neumann and Larry Siebenmann:

 \catcode`\@=11
 \let\wlog@ld\wlog
 \def\wlog#1{\relax}

 \newif\ifIN@
 \def\m@rker{\m@@rker}
 \def\IN@{\expandafter\INN@\expandafter}
 \long\def\INN@0#1@#2@{\long\def\NI@##1#1##2##3\ENDNI@
    {\ifx\m@rker##2\IN@false\else\IN@true\fi}%
     \expandafter\NI@#2@@#1\m@rker\ENDNI@}
  \newtoks\Initialtoks@  \newtoks\Terminaltoks@
  \def\SPLIT@{\expandafter\SPLITT@\expandafter}
  \def\SPLITT@0#1@#2@{\def\TTILPS@##1#1##2@{%
     \Initialtoks@{##1}\Terminaltoks@{##2}}\expandafter\TTILPS@#2@}
  \newtoks\Trimtoks@

 \def\ForeTrim@{\expandafter\ForeTrim@@\expandafter}
 \def\ForePrim@0 #1@{\Trimtoks@{#1}}
 \def\ForeTrim@@0#1@{\IN@0\m@rker. @\m@rker.#1@%
     \ifIN@\ForePrim@0#1@%
     \else\Trimtoks@\expandafter{#1}\fi}
 
  \def\Trim@0#1@{%
      \ForeTrim@0#1@%
      \IN@0 @\the\Trimtoks@ @%
        \ifIN@
             \SPLIT@0 @\the\Trimtoks@ @\Trimtoks@\Initialtoks@
             \IN@0\the\Terminaltoks@ @ @%
                 \ifIN@
                 \else \Trimtoks@ {FigNameWithSpace}%
                 \fi
        \fi
      }

  %%% Math Bolds
  \font\titlebold=cmbx12 scaled 1200
  \font\twelvebold=cmbx12
  \font\tenbold=cmbx10
  \font\ninebold=cmbx9
  \font\sevenbold=cmbx7
  \font\fivebold=cmbx5

  \input amssym.def \input amssym
  %%% point sizes not loaded by amssym.def:
     \font\titlemsa=msam10 at 14.4pt
     \font\titlemsb=msbm10 at 14.4pt
     \font\titleeufm=eufm10 at 14.4pt
     \font\twelvemsa=msam10 scaled 1200
     \font\twelvemsb=msbm10 scaled 1200
     \font\twelveeufm=eufm10 scaled 1200
     \font\ninemsa=msam9
     \font\ninemsb=msbm9
     \font\nineeufm=eufm9

   %%% Cyrillic fonts (for accents and input, see ams cyr doc)
   \ifx\cyrfam\undefined
   \else
     \immediate\write16{}%
     \message{ !!! cyr fonts already defined. !!! }
     \message{ --- edit out superfluous font defs? }
   \fi
   \newfam\cyrfam
       \font\titlecyr=wncyr10 scaled 1440 %%% no caps?
       \font\twelvecyr=wncyr10 scaled 1200
       \font\tencyr=wncyr10
       \font\ninecyr=wncyr9
       \font\sevencyr=wncyr7
       \font\sixcyr=wncyr6

   %%% Euler script fonts (replacing caligraphic):
   \newfam\eusmfam
       \font\titleeusm=eusm10 scaled 1440
       \font\twelveeusm=eusm10 scaled 1200
       \font\teneusm=eusm10
       \font\nineeusm=eusm9
       \font\seveneusm=eusm7
       
       \font\fiveeusm=eusm5

\let\Cal\cal

 %%% Some fonts not loaded by plain
    \font\ninemrm=cmr9 %% new name for 9 pt math roman
    \font\ninei=cmmi9
    \font\ninesy=cmsy9 
    \skewchar\ninei='177
    \skewchar\ninesy='60

  \font\twelvemrm=cmr10 at 12pt %% new name
  \font\twelvei=cmmi10 at 12pt
  \font\twelvesy=cmsy10 at 12pt
 % \font\twelveex=cmex10 at 12pt

  \font\titlemrm=cmr10 at 14.4pt %% new name
  \font\titlei=cmmi10 at 14.4pt
  \font\titlesy=cmsy10 at 14.4pt
 % \font\titleex=cmex10 at 14.4pt

 %%%% Miscellanious font definitions

  \def\Smallfonts{\ninepoint}

  \def\Hfont{\titlepoint\bf}
  \def\Authorfont{\twelvepoint\it}
  \def\HHfont{\twelvepoint\bf}
  \def\HHHfont{\bf}
  % automatically smaller in 9 point parts
  \def\Bibfont{\tenbf}
  \def\Coordfont{\nineit }% defined in osuPSfnt.sty

  \def \thfont {\bf }
  \def \pffont {\it\itSpacing }
  \def \rkfont {\bf }
  \def \dffont {\bf }
  \def \egfont {\bf }

 %%%%% NINEPOINT %%%%%
 \def\ninepoint{%
  \def\rm{\fam0\ninerm}%
    \textfont0=\ninemrm  \scriptfont0=\sevenrm  \scriptscriptfont0=\fiverm
    \textfont1=\ninei    \scriptfont1=\seveni   \scriptscriptfont1=\fivei
  \def\mit{\fam1\ninei}%
  \def\oldstyle{\fam1\ninei}%
    \textfont2=\ninesy   \scriptfont2=\sevensy  \scriptscriptfont2=\fivesy
    \textfont3=\tenex    \scriptfont3=\tenex    \scriptscriptfont3=\tenex
  \def\it{\fam\itfam\nineit}%
    \textfont\itfam=\nineit
  \def\bf{\ifmmode\fam\bffam\else\ninebf\fi}%
    \textfont\bffam=\ninebold 
    \scriptfont\bffam=\sevenbold 
    \scriptscriptfont\bffam=\fivebold%
  \def\msa{\fam\msafam\ninemsa}%
    \textfont\msafam=\ninemsa 
    \scriptfont\msafam=\sevenmsa
    \scriptscriptfont\msafam=\fivemsa%
  \def\msb{\fam\msbfam\ninemsb}%
    \textfont\msbfam=\ninemsb%
    \scriptfont\msbfam=\sevenmsb%
    \scriptscriptfont\msbfam=\fivemsb%
  \def\eufm{\fam\eufmfam\nineeufm}%
    \textfont\eufmfam=\nineeufm
    \scriptfont\eufmfam=\seveneufm
    \scriptscriptfont\eufmfam=\fiveeufm
   \def\eusm{\fam\eusmfam\nineeusm}%
     \textfont\eusmfam=\nineeusm
     \scriptfont\eusmfam=\seveneusm
     \scriptscriptfont\eusmfam=\fiveeusm
   \def\cyr{\fam\cyrfam\ninecyr}%
     \textfont\cyrfam=\ninecyr
     \scriptfont\cyrfam=\sevencyr
     \scriptscriptfont\cyrfam=\sixcyr%%
  \setbox\strutbox=\hbox{\vrule
      height7pt depth3pt width0pt}%
   \baselineskip=10.8pt\rm}

 \let\eightpoint\ninepoint % we do not use eightpoint

 %%%%% FONTS AT TENPOINT %%%%%
 \def\tenpoint{%
  \def\rm{\fam0\tenrm}%
    \textfont0=\tenmrm \scriptfont0=\sevenrm \scriptscriptfont0=\fiverm%
  \def\mit{\fam1\teni}%
  \def\oldstyle{\fam1\teni}%
    \textfont1=\teni   \scriptfont1=\seveni  \scriptscriptfont1=\fivei%
    \textfont2=\tensy  \scriptfont2=\sevensy \scriptscriptfont2=\fivesy%
    \textfont3=\tenex  \scriptfont3=\tenex   \scriptscriptfont3=\tenex%
  \def\it{\fam\itfam\tenit}%
    \textfont\itfam=\tenit%
  \def\bf{\ifmmode\fam\bffam\else\tenbf\fi}%
    \textfont\bffam=\tenbold% was tenbold for osu
    \scriptfont\bffam=\sevenbold%
    \scriptscriptfont\bffam=\fivebold%
  \def\msa{\fam\msafam\tenmsa}%
    \textfont\msafam=\tenmsa%
    \scriptfont\msafam=\sevenmsa%
    \scriptscriptfont\msafam=\fivemsa%
  \def\msb{\fam\msbfam\tenmsb}%
    \textfont\msbfam=\tenmsb%
    \scriptfont\msbfam=\sevenmsb%
    \scriptscriptfont\msbfam=\fivemsb%
  \def\eufm{\fam\eufmfam\teneufm}%
   \textfont\eufmfam=\teneufm
   \scriptfont\eufmfam=\seveneufm
   \scriptscriptfont\eufmfam=\fiveeufm
   \def\eusm{\fam\eusmfam\teneusm}%
    \textfont\eusmfam=\teneusm
    \scriptfont\eusmfam=\seveneusm
    \scriptscriptfont\eusmfam=\fiveeusm
   \def\cyr{\fam\cyrfam\tencyr}%
    \textfont\cyrfam=\tencyr
    \scriptfont\cyrfam=\sevencyr
    \scriptscriptfont\cyrfam=\sixcyr%%
  \setbox\strutbox=\hbox{\vrule %
      height8.5pt depth3.5ptwidth0pt}%
  \baselineskip=\StdBaselineskip\rm}

 %%%%% FONTS AT TWELVEPOINT %%%%%
 \def\twelvepoint{%
  \def\rm{\fam0\twelverm}%
    \textfont0=\twelvemrm \scriptfont0=\tenmrm \scriptscriptfont0=\sevenrm
    \textfont1=\twelvei   \scriptfont1=\teni   \scriptscriptfont1=\seveni
  \def\mit{\fam1\twelvei}%
  \def\oldstyle{\fam1\twelvei}%
    \textfont2=\twelvesy  \scriptfont2=\tensy  \scriptscriptfont2=\sevensy
    \textfont3=\tenex  \scriptfont3=\tenex  \scriptscriptfont3=\tenex
  \def\it{\fam\itfam\twelveit}%
    \textfont\itfam=\twelveit
  \def\bf{\ifmmode\fam\bffam\else\twelvebf\fi}%
    \textfont\bffam=\twelvebold
    \scriptfont\bffam=\tenbold%
    \scriptscriptfont\bffam=\sevenbold%
  \def\msa{\fam\msafam\twelvemsa}%
    \textfont\msafam=\twelvemsa%
    \scriptfont\msafam=\tenmsa%
    \scriptscriptfont\msafam=\sevenmsa%
  \def\msb{\fam\msbfam\twelvemsb}%
    \textfont\msbfam=\twelvemsb%
    \scriptfont\msbfam=\tenmsb%
    \scriptscriptfont\msbfam=\sevenmsb%
  \def\eufm{\fam\eufmfam\twelveeufm}%
   \textfont\eufmfam=\twelveeufm
   \scriptfont\eufmfam=\teneufm
   \scriptscriptfont\eufmfam=\seveneufm
   \def\eusm{\fam\eusmfam\twelveeusm}%
    \textfont\eusmfam=\twelveeusm
    \scriptfont\eusmfam=\teneusm
    \scriptscriptfont\eusmfam=\seveneusm
   \def\cyr{\fam\cyrfam\tencyr}%
    \textfont\cyrfam=\twelvecyr
    \scriptfont\cyrfam=\tencyr
    \scriptscriptfont\cyrfam=\sevencyr%%
  \setbox\strutbox=\hbox{\vrule
      height10.2pt depth4.55pt width0pt}%
  \baselineskip=14pt\rm}

 %%%%% FONTS AT TITLEPOINT %%%%%
 \def\titlepoint{%
    \textfont0=\titlemrm \scriptfont0=\twelvemrm \scriptscriptfont0=\tenmrm
    \textfont1=\titlei   \scriptfont1=\twelvei   \scriptscriptfont1=\teni
  \def\mit{\fam1\titlei}%
  \def\oldstyle{\fam1\titlei}%
    \textfont2=\titlesy  \scriptfont2=\twelvesy  \scriptscriptfont2=\tensy
    \textfont3=\tenex% math ext not avail in varying sizes??
    \scriptfont3=\tenex
    \scriptscriptfont3=\tenex
  \def\it{\fam\itfam\titleit}%
    \textfont\itfam=\titleit
  \def\bf{\ifmmode\fam\bffam\else\titlebf\fi}%
    \textfont\bffam=\titlebold
    \scriptfont\bffam=\twelvebold%
    \scriptscriptfont\bffam=\tenbold%
  \def\msa{\fam\msafam\titlemsa}%
    \textfont\msafam=\titlemsa%
    \scriptfont\msafam=\twelvemsa%
    \scriptscriptfont\msafam=\tenmsa%
  \def\msb{\fam\msbfam\titlemsb}%
    \textfont\msbfam=\titlemsb%
    \scriptfont\msbfam=\twelvemsb%
    \scriptscriptfont\msbfam=\tenmsb%
  \def\eufm{\fam\eufmfam\titleeufm}%
    \textfont\eufmfam=\titleeufm
    \scriptfont\eufmfam=\twelveeufm
    \scriptscriptfont\eufmfam=\teneufm
   \def\eusm{\fam\eusmfam\titleeusm}%
     \textfont\eusmfam=\titleeusm
     \scriptfont\eusmfam=\twelveeusm
     \scriptscriptfont\eusmfam=\teneusm
   \def\cyr{\fam\cyrfam\tencyr}%
    \textfont\cyrfam=\titlecyr
    \scriptfont\cyrfam=\twelvecyr
    \scriptscriptfont\cyrfam=\tencyr%%
  \setbox\strutbox=\hbox{\vrule
      height12.3pt depth5.54pt width0pt}%
  \baselineskip=16pt\rm}

 %%%% RUNNING HEADINGS
\newbox\AuthorBox\newbox\TitleBox
\newbox\TFLinebox
\newbox\FLinebox
\newbox\HLinebox
\def\SetTFLinebox#1{\setbox\TFLinebox=\hbox{#1}}
\def\SetFLinebox#1{\setbox\FLinebox=\hbox{#1}}
\def\SetHLinebox#1{\setbox\HLinebox=\hbox{#1}}

 \def\SetAuthorHead#1{%
     \setbox\AuthorBox=\hbox{\ninepoint \it 
           \ignorespaces\frenchspacing#1\unskip}}
 \def\SetTitleHead#1{%
     \setbox\TitleBox=\hbox{\ninepoint \it
           \ignorespaces\frenchspacing#1\unskip}}

 %% Italic Spacing Correction
  \def\itSpacing{\relax}
  \def\itSpacingOff{\relax}

  %% Main section headings

 \def\Hrule{\hrule width0pt height0pt}

 %% skip used around proclamations, after section headings,
  % and before subsection-headings:
  \newskip\ProcSkip \ProcSkip 8pt plus2pt minus2pt

 \newskip\LastSkip
 \def\SaveLastSkip{\LastSkip\lastskip}
 \def\RestoreLastSkip{\vskip-\LastSkip\vskip\LastSkip}

 %% Do not indent next paragraph after a header:
 \def\NoindentAfter{\everypar={\setbox0=\lastbox\everypar={}}}

 \long\def\H#1\par#2\par{\notenumber=0 \titlepagetrue%
    {
    \baselineskip=20pt
    \parindent=0pt\parskip=0pt\frenchspacing
    \leftskip=0pt plus .2\hsize minus .3\hsize
    \rightskip=0pt plus .2\hsize minus .3\hsize
 \def\\{\unskip\break}%
    \pretolerance=10000 \Hfont #1\unskip\break
     \vskip7pt\Hrule
\hfill \Authorfont #2\hfill\hfill\unskip}
    \vskip48pt plus 4pt minus 4pt% 60pt=48+12pt
    \par\NoindentAfter\rm}

 \long\def\Hi#1\par#2\par{\notenumber=0 \titlepagetrue%
    {  \baselineskip=0pt  \parindent=0pt\parskip=0pt\frenchspacing
    \leftskip=0pt plus .2\hsize minus .3\hsize
    \rightskip=0pt plus .2\hsize minus .3\hsize
}
    \rm}

 %%% Minor section headings

 \newdimen\PageRemainder
  \def\SetPageRemainder{%\maxdimen case at page tops 12-91 LS
     \PageRemainder=\pagegoal
     \ifdim\PageRemainder=\maxdimen\PageRemainder=\vsize
     \else\advance\PageRemainder by -1\pagetotal\fi}

  \def\Rpt@{}\def\Rpt@@{}

  \long\def\HH#1\par{\par%A
  \SaveLastSkip\removelastskip\goodbreak
  \ifdim\LastSkip<30pt %24pt
     \LastSkip 30pt%24pt 
plus 3pt minus 2pt\fi
  \SetPageRemainder\advance\PageRemainder-\LastSkip
  \ifdim\PageRemainder<150pt
       \edef\Rpt@{remain = \the\PageRemainder\noexpand\\
                pagetotal=\the\pagetotal\noexpand\\
                           pagegoal=\the\pagegoal}%
          \fi
   \ifdim\PageRemainder<65pt %%Head plus 4 lines (approx)
       \ifdim\PageRemainder > 0pt
          \edef\Rpt@@{\noexpand\\
                      Had HH PageRemainder$<$\relax 65pt\noexpand\\
                      Hence forced break!}%
     \vskip 0pt plus .2\PageRemainder\eject %% Pull it out a bit
    \fi\fi
    \vskip\LastSkip\Hrule %%%%%%%%\Hrule added
    \pretolerance=10000\rightskip=0pt plus 3em%B
    \hangafter1 \hangindent=2.2em%
    \noindent
    \HHfont \unskip \Ednote{\Rpt@\Rpt@@}%
            \def\Rpt@{}\def\Rpt@@{}%
            \ignorespaces
            #1\par\rightskip=0pt\pretolerance=\StdPretolerance%
    \NoindentAfter
\tenpoint\rm%
     \medskip \vskip\ProcSkip}%interlineskip adds 2pt to this

  \long\def\HHH#1\par{\par%
  \SaveLastSkip\removelastskip\goodbreak
  \ifdim\LastSkip<\ProcSkip%
     \LastSkip\ProcSkip\fi
  \SetPageRemainder\advance\PageRemainder-\LastSkip
  \ifdim\PageRemainder<150pt
       \edef\Rpt@{remain = \the\PageRemainder\noexpand\\
                pagetotal=\the\pagetotal\noexpand\\
                           pagegoal=\the\pagegoal}%
       \fi
   \ifdim\PageRemainder<48pt  %% 4 lines
        \ifdim\PageRemainder > 0pt
             \edef\Rpt@@{\noexpand\\
                      Had HHH PageRemainder$<$\relax48pt\noexpand\\
                      Hence forced break!}%
       \vskip 0pt plus .2\PageRemainder\eject %% Pull it out a bit
      \fi\fi
   \vskip\LastSkip\par\noindent
   \HHHfont \unskip\Ednote{\Rpt@\Rpt@@}%
  \def\Rpt@{}\def\Rpt@@{}%
  \ignorespaces
   #1\unskip.\quad\rm\ignorespaces
   \ignorepars}

  \long\def\ignorepars#1\par{\def\Test{#1}%
     \ifx\Test\Empty\def\This{\ignorepars}%
        \else\def\This{\Test\par}\fi
           \This}
  \def\Empty{}

 \def\Abstract#1\par{\bgroup\Smallfonts\narrower\HHH #1\par}
 \def\endAbstract{\par\egroup}

 %%%%% Proclamations %%%%%

 \def\ProcBreak{\par%
    \ifdim\lastskip<8pt%
    \removelastskip%
    \penalty-200\vskip\ProcSkip\fi}

 \def\th#1\par{\ProcBreak \noindent
   {\thfont\ignorespaces
    #1\unskip.}\it\itSpacing\kern.4em\ignorepars}%\everymath{\/}

 \def\endth{\ProcBreak\rm\itSpacingOff }%\everymath{}

  %% the theorem statement will be in italic by default

 \def\pf#1\par{\ProcBreak %
    \noindent\pffont#1\unskip.\rm\itSpacingOff{\kern .7em}\ignorepars}

  %% \qed is alternative

  %% A Box for the QED
  \def\qedbox{\hbox{\vbox{
    \hrule width0.2cm height0.2pt
    \hbox to 0.2cm{\vrule height 0.2cm width 0.2pt
             \hfil\vrule height0.2cm width 0.2pt}
    \hrule width0.2cm height 0.2pt}\kern1pt}}

  %% Typing in \qed makes the qedbox right justified:
  \def\qed{\ifmmode\qedbox
    \else\unskip\ \hglue0mm\hfill\qedbox\ProcBreak\fi}

  \def \rk #1\par{\ProcBreak
     \noindent{\rkfont\ignorespaces #1\unskip.}%
     \rm\kern.6em\ignorepars}

  \def \df #1\par{\ProcBreak
     \noindent{\dffont\unskip\ignorespaces #1\unskip.}%
     \rm\kern.6em\ignorepars}

  \def \enddf {\medskip\ProcBreak }

  \def \eg #1\par{\ProcBreak
     \noindent\egfont\unskip\ignorespaces #1\unskip.
     \rm\kern.6em\ignorepars}

  \def \endeg {\medskip\ProcBreak }

  \newdimen\Overhang

   \def\MaxTag@#1#2#3#4#5{\setbox0=\hbox{#4\ignorespaces#2\unskip}%
     \dimen0=\wd0\advance\dimen0 by#3
     \ifdim\dimen0<#5\relax\dimen0=#5\fi
     \expandafter\edef\csname #1Hang\endcsname{\the\dimen0}}

 \def\MaxItemTag#1{\MaxTag@{Item}{#1}{.4em}{\ItemStyle}{\parindent}}%
 \def\MaxItemItemTag#1{%
        \MaxTag@{ItemItem}{#1}{.4em}{\ItemItemStyle}{\parindent}}
 \def\MaxNrTag#1{\MaxTag@{Nr}{#1}{.5em}{\NrStyle}{\parindent}}
 \def\MaxReferenceTag#1{%
        \MaxTag@{Reference}{[#1]}{.6em}{\ninerm}{\parindent}}
 \def\MaxFootTag#1{\MaxTag@{Foot}{#1}{.4em}{\ninerm}{\z@}}

  %% \SetOverhang@ will prevent for tag-text collision
  \def\SetOverhang@{\Overhang=.8\dimen0%
     \advance\Overhang by \wd0\relax%nec!
     \ifdim\Overhang>\hangindent\relax%nec!
       \advance\Overhang by .25\dimen0%
       \Ednote{Tag is pushing text.}\osumess{Tag is pushing text.}%
     \else\Overhang=\hangindent
     \fi}

   %%% \Item
   \def\Item#1{\par\noindent
      \hangafter1\hangindent=\ItemHang
      \setbox0=\hbox{\ItemStyle\ignorespaces#1\unskip}%
      \dimen0=.4em\SetOverhang@% dimen0 is extra space
      \rlap{\box0}\kern\Overhang\ignorespaces}

   %%% \ItemItem
   \def\ItemItem#1{\par\noindent
      \hangafter1\hangindent=\ItemItemHang
      \setbox0=\hbox{\ItemItemStyle\ignorespaces#1\unskip}%
      \dimen0=.4em\SetOverhang@
      \advance\hangindent by \ItemHang
      \kern\ItemHang\rlap{\box0}%
      \kern\Overhang\ignorespaces}

  %%%% \Nr Items without hanging indentation
  \def\Nr#1{\par\noindent\hangindent=\NrHang %not really a hang
    \setbox0=\hbox{\NrStyle\ignorespaces#1\unskip}%
    \dimen0=.5em\SetOverhang@% dimen0 is extra space
    \rlap{\box0}\kern\Overhang
    \hangindent=\z@\ignorespaces}

  %%%% Roster (not compulsory)
  %%  endRoster has to remember \lastskip (e.g. from a \qed) through \egroup.
   \newskip\Rosterskip\Rosterskip 1pt plus1pt %% modifiable
   \def\Roster{\par\ifdim\lastskip<\Rosterskip\removelastskip\vskip\Rosterskip\fi
    \bgroup}
   \def\endRoster{\par\global\edef\LastSkip@{\the\lastskip}\removelastskip
       \egroup\penalty-50\LastSkip\LastSkip@\relax
       \ifdim\LastSkip<\Rosterskip\LastSkip\Rosterskip\fi
       \vskip\LastSkip}%%changed Feb/5/92 WN

 %%%%% Emphasis %%%%%

 %%%%% Vertical spacing %%%%%

 %%%%% References %%%%%

 \def\cite#1{%\relaxnext@
    \def\nextiii@##1,##2\end@{{\frenchspacing\rm 
      \lBr\ignorespaces##1\unskip{\rm,~\ignorespaces##2}\rBr}}%
    \IN@0,@#1@%
    \ifIN@\def\next{\nextiii@#1\end@}\else
    \def\next{{\rm\lBr#1\rBr}}\fi\next}

 %%%%% Bibliography %%%%%

   \def \Bib#1\par{%
       \par\removelastskip\SetPageRemainder
       \ifdim\PageRemainder < 97pt
        \ifdim\PageRemainder > 0pt
        \vfill\eject
       \fi\fi
    \ProcBreak \par\begingroup\parskip=0 pt%
    \goodbreak \vskip 15 pt plus 10 pt
    \noindent\null\hfill\Bibfont% \kern??pt%  (center over what?)
      \ignorespaces #1\unskip\hfill\null\par 
    \frenchspacing \Smallfonts\rm
    \parskip=2.5 pt plus 1 pt minus.5pt%
    \nobreak\vskip 12pt plus 2pt minus2pt\nobreak
    \leftskip=0 pt \baselineskip=10.5pt}

 \def\ReferenceTagSlide{0em}
  \def\ReferenceTagGap{.5em}

  \def \rf#1{\par\noindent
     \hangafter1\hangindent=\ReferenceHang      
     \setbox0=\hbox{\ninerm[\ignorespaces#1\unskip]}%        
     \dimen0=\ReferenceTagGap\SetOverhang@
     \rlap{\kern\ReferenceTagSlide\box0}%       
     \kern\Overhang\ignorespaces}

  \def\ref#1\par#2\par#3\par#4\par{%
     \rf{#1}#2\unskip,\ #3\unskip,\
     #4\unskip.}

  \def\endBib{\par\endgroup\vskip 12pt minus 6pt }

 %%%%% Coordinates %%%%%

  \long\def\Coordinates#1\endCoordinates{%\relax}
 {\par\vskip4pt\def\\{\unskip, }\Coordfont\baselineskip10.5pt\noindent#1}}

 \def\pagecontents{%\TRMargIns new, \Pagetot@l
  \gdef\Pagetot@l{\pagetotal}
  \ifvoid\TRMargIns\else
    \rlap{\kern\hsize\kern10pt\vbox to 0pt{%
         \box\TRMargIns\vss}}\fi
  \ifvoid\topins\else\unvbox\topins\fi
   \dimen@=\dp\@cclv \unvbox\@cclv % open up \box255
   \ifvoid\footins\else % footnote info is present
     \vskip\skip\footins
     \footnoterule
     \unvbox\footins\fi
   \ifr@ggedbottom \kern-\dimen@ \vfil \fi}

  %%%%% Some math accents %%%%%

 \newcount\Ht %pg121; Height register, used in Linefigure & accents

 \def \Acc{\expandafter } %%% What is this for?? WN

 \def\swthat{\raise -1.1 ex\hbox{\sam$\widehat{}$}}
 \def\swttilde{\raise -1.2 ex\hbox{\sam$\widetilde{}$}}
 \def \overdot{{\raise .2 ex \hbox to 0pt {\hss\bf\smash{.}\hss}}}
 \def \overcircle{{\raise .1 ex \hbox to 0pt
    {\sam$\eightpoint\scriptstyle\hss\circ\hss$}}}

 \def \Mathaccent#1#2{{\sam % E.g. #1=\widehat
  \setbox4=\hbox{$\vphantom{#2}$}
  \Ht=\ht4 %pg120
  \setbox5=\hbox{${#1}$}
  \setbox6=\hbox{${#2}$}
  \setbox7=\hbox to .5\wd6{}
  \copy7\kern .1\Ht \raise\Ht sp\hbox{\copy5}\kern-.1\Ht
  \copy7\llap{\box6}
  }}

  \def\SwtCheck #1{
        \ifmmode \check{#1}%
                \else \v {#1}%
                \fi}

 %%  \barpartial : bar over partial is common, tailor!
 \def\barpartial {%
   \kern .17 em
    \overline {\kern -.17 em\partial\kern-.03 em}%
    \kern .03 em}

 %%%   BEtter overline
 
  \def\Overline#1{\setbox1=\hbox{\sam ${#1}$}%
      \ifdim \wd1 > 6pt
    \kern .11 em
    \overline {\kern -.11 em#1\kern-.14 em}
    \kern .14 em
  \else
    \kern .03 em
    \overline {\kern -.03 em#1\kern-.04 em}
    \kern .04 em
  \fi}

 \def\SOverline#1{\setbox1=\hbox{\sam ${#1}$}%
      \ifdim \wd1 > 7pt
    \kern .22 em
    \overline {\kern -.22 em#1\kern-.09 em}%
    \kern .09 em
  \else
    \kern .10 em
    \overline {\kern -.10 em#1\kern-.04 em}%
    \kern .04 em
  \fi}

  %%% Better underline

 \def\Underline#1{\setbox1=\hbox{\sam ${#1}$}%
      \ifdim \wd1 > 6pt
    \kern .11 em
    \underline {\kern -.11 em#1\kern-.14 em}
    \kern .14 em
  \else
    \kern .03 em
    \underline {\kern -.03 em#1\kern-.04 em}
    \kern .04 em
  \fi}

 \def\SUnderline#1{\setbox1=\hbox{\sam ${#1}$}%
      \ifdim \wd1 > 7pt
    \kern .04 em
    \underline {\kern -.04 em#1\kern-.2 em}%
    \kern .2 em
  \else
    \kern .0 em
    \underline {\kern -.0 em#1\kern-.15 em}%
    \kern .15 em
  \fi}

  %%%%% Miscellaneous %%%%%

 \def \Blackbox
   {\leavevmode\hskip .3pt \vbox
   {\hrule height 5pt\hbox{\hskip 4.5pt}}\hskip .5pt}

 \def \XX{\Blackbox\kern.5pt\Blackbox} %% editorial use

  \def\.{.\kern1pt}

  %% unbreakable hyphen (by local change of hyphenchar to -1)
    \def\Hyphen{\edef\this{\the\hyphenchar\font}%
          \hyphenchar\font=-1\char\this\hyphenchar\font=\this}

  %% Prose In Math or Display 
 \ifx\undefined\text
  \def\text#1{\hbox{\rm #1}}\fi %% AMSTeX is more sophisticated

  %% Math Object Names (multi-character math object names)
  %%\nolimits can be cancelled
                                     % by a following \limits if wanted

%%%% Larry's mathsurround hacks:

   \everymath{}  %% initially, but later ...

  \def\PassMath@@{\aftergroup\AfterMath@} %% use \aftergroup LS 5-92

 \let\PassMath@\PassMath@@

 \def\AfterMath@{\futurelet\next\AfterMathMole@}

 \def\AfterMathMole@{%\show\next
      \ifcat\next\space% picks off CR and \par cases too; not \dots
          \def\this{}%{(space)}%
      \else
      \ifcat\next\egroup %
        \def\this{\osumess{Handset mathsurround?? ---(see dollar brace)}}%
      \else
      \def\this{\AAfterMath@}% this minority case slow
      \fi\fi
      \this}

 \def\hyphen@{-}
 \def\paren@{)}
 \def\apostr@{'}

 \def\MSC#1{\kern-.8\mathsurround#1\kern.8\mathsurround}

 \def\AAfterMath@#1{\def\Next{#1}%\show\Next%
    \IN@0\Next @,.;:!?\relax @%
    \ifIN@\def\this{\MSC{\Next}}%
    \else
    \ifx\Next\hyphen@\def\this{\futurelet\next\AfterHyphen@}%
    \else
    \ifx\Next\paren@\def\this{#1}%
    \else 
    \ifx\Next\apostr@\def\this{#1}%
    \else \def\this{\osumess{Handset mathsurround??}%
                 #1}\fi\fi\fi\fi
    \this}

 \def\AfterHyphen@#1{\def\Next{#1}%
   \ifx\Next\hyphen@\def\this{--}\else
   \ifcat\next\space%
   \def\this{\kern-\mathsurround\kern.05em- \Next}\else
   \def\this{\kern-\mathsurround\kern.05em\Hyphen\Next}\fi\fi\this}

%%%% switches
 \def\sam{\mathsurround=\z@\let\PassMath@\relax}  %
 \def\mas{\mathsurround=\StdMathsurround\let\PassMath@\PassMath@@}
 
 \def\Mas{\mathsurround=\StdMathsurround
                \everymath{\PassMath@}\let\PassMath@\PassMath@@}

 \def\m@th{\mathsurround=\z@\everymath{}}%% good general measure

 \def\m@@th{\mathsurround=\z@\everymath={}\let\m@th\relax}

\def\underbar#1{$\setbox\z@\hbox{#1}\dp\z@\z@
      \m@th \underline{\box\z@}$\relax}

\def\mathhexbox#1#2#3{\leavevmode
  \hbox{\m@@th$\m@th \mathchar"#1#2#3$}}

\def\dots{\relax\ifmmode\ldots\else$\m@th\ldots\,$\relax\fi}
   %%% this first \relax is ONLY original

\def\dotfill{\cleaders\hbox{\m@@th$\m@th \mkern1.5mu.\mkern1.5mu$}\hfill}
\def\rightarrowfill{$\m@th\mathord-\mkern-6mu%
  \cleaders\hbox{\m@@th$\mkern-2mu\mathord-\mkern-2mu$}\hfill
  \mkern-6mu\mathord\rightarrow$\relax}
\def\leftarrowfill{$\m@th\mathord\leftarrow\mkern-6mu%
  \cleaders\hbox{\m@@th$\mkern-2mu\mathord-\mkern-2mu$}\hfill
  \mkern-6mu\mathord-$\relax}

\def\downbracefill{$\m@th\braceld\leaders\vrule\hfill\braceru
  \bracelu\leaders\vrule\hfill\bracerd$\relax}
\def\upbracefill{$\m@th\bracelu\leaders\vrule\hfill\bracerd
  \braceld\leaders\vrule\hfill\braceru$\relax}

\def\angle{{\vbox{\m@@th\ialign{$\m@th\scriptstyle##$\crcr
      \not\mathrel{\mkern14mu}\crcr
      \noalign{\nointerlineskip}
      \mkern2.5mu\leaders\hrule height.34pt\hfill\mkern2.5mu\crcr}}}}

\def\big#1{{\m@@th\hbox{$\left#1\vbox to8.5\p@{}\right.\n@space$}}}
\def\Big#1{{\m@@th\hbox{$\left#1\vbox to11.5\p@{}\right.\n@space$}}}
\def\bigg#1{{\m@@th\hbox{$\left#1\vbox to14.5\p@{}\right.\n@space$}}}
\def\Bigg#1{{\m@@th\hbox{$\left#1\vbox to17.5\p@{}\right.\n@space$}}}
\def\n@space{\nulldelimiterspace\z@ \m@th}

\def\root#1\of{\setbox\rootbox\hbox{\m@@th$\m@th\scriptscriptstyle{#1}$}
  \mathpalette\r@@t}
\def\r@@t#1#2{\setbox\z@\hbox{\m@@th$\m@th#1\sqrt{#2}$\relax}
  \dimen@\ht\z@ \advance\dimen@-\dp\z@
  \mkern5mu\raise.6\dimen@\copy\rootbox \mkern-10mu \box\z@}

\def\mathph@nt#1#2{\setbox\z@\hbox{\m@@th$\m@th#1{#2}$}\finph@nt}

\def\mathsm@sh#1#2{\setbox\z@\hbox{\m@@th$\m@th#1{#2}$}\finsm@sh}

\def\@vereq#1#2{\lower.5\p@\vbox{\m@@th\baselineskip\z@skip\lineskip-.5\p@
    \ialign{$\m@th#1\hfil##\hfil$\crcr#2\crcr=\crcr}}}

\def\mathpalette#1#2{\sam\mathchoice{#1\displaystyle{#2}}%
  {#1\textstyle{#2}}{#1\scriptstyle{#2}}{#1\scriptscriptstyle{#2}}\mas}

\def\widehat#1{\setbox\z@\hbox{\sam$#1$}%
 \ifdim\wd\z@>\tw@ em\mathaccent"0\msbfam@5B{#1}%
 \else\mathaccent"0362{#1}\fi}
\def\widetilde#1{\setbox\z@\hbox{\sam$#1$}%
 \ifdim\wd\z@>\tw@ em\mathaccent"0\msbfam@5D{#1}%
 \else\mathaccent"0365{#1}\fi}

 \def\dots{\relax{}
  \ifmmode\def\thedots{\mdots@}\else\def\thedots{\tdots@}\fi %
  \thedots}

 %% \eqno and \leqno need protection
 \let\@ldeqno\eqno\let\@ldleqno\leqno
 \def\eqno{\everymath{}\@ldeqno} \def\leqno{\everymath{}\@ldleqno}

  \let\@ldeqalignno\eqalignno
  \def\eqalignno#1{\sam\@ldeqalignno{#1}\mas}
  \let\@ldeqalign\eqalign
  \def\eqalign#1{\sam\@ldeqalign{#1}\mas}

 \def\overrightarrow#1{\vbox{\m@th\ialign{##\crcr
      \rightarrowfill\crcr\noalign{\kern-\p@\nointerlineskip}
      $\hfil\displaystyle{#1}\hfil$\crcr}}}
 \def\overleftarrow#1{\vbox{\m@th\ialign{##\crcr
      \leftarrowfill\crcr\noalign{\kern-\p@\nointerlineskip}
      $\hfil\displaystyle{#1}\hfil$\crcr}}}
 \def\overbrace#1{\mathop{\vbox{\m@th\ialign{##\crcr\noalign{\kern3\p@}
      \downbracefill\crcr\noalign{\kern3\p@\nointerlineskip}
      $\hfil\displaystyle{#1}\hfil$\crcr}}}\limits}
 \def\underbrace#1{\mathop{\vtop{\m@th\ialign{##\crcr
      $\hfil\displaystyle{#1}\hfil$\crcr\noalign{\kern3\p@\nointerlineskip}
      \upbracefill\crcr\noalign{\kern3\p@}}}}\limits}

  \let\@ldmatrix\matrix
  \let\end@ldmatrix\endmatrix
  \def\matrix{\sam\@ldmatrix}
  \def\endmatrix{\end@ldmatrix\mas}
  \let\@ldgather\gather
  \let\end@ldgather\endgather
  \def\gather{\sam\@ldgather}
  \def\endgather{\end@ldgather\mas}
  \let\@ldalign\align
  \let\end@ldalign\endalign
  \def\align{\sam\@ldalign}
  \def\endalign{\end@ldalign\mas}
  \let\@ldaligned\aligned
  \let\end@ldaligned\endaligned
  \def\aligned{\sam\@ldaligned}
  \def\endaligned{\end@ldaligned\mas}
  \let\@ldtag\tag
  \def\tag{\sam\@ldtag}
   %
  %%% Commutative diagrams : use LamsCD too?

   \let\MinCDArrowWidth\minCDaw@

  %% will be redefined by BoxedEPS.tex

  %%%%% \FigureTitle %%%%%

%%%% End of Larry's mathsurround stuff
%%%% Start of Walter's insert corrections

\newskip\insertskipamount\newskip\inserthardskipamount
\insertskipamount 6pt plus2pt %This is medskipamount without shrink
\inserthardskipamount 6pt
\def\insertskip{\vskip\insertskipamount}
\newcount\SplitTest%        will be set to -1 if a topinsert has split
\def\SetSplitTest{\SplitTest\insertpenalties
  \insert\topins{\floatingpenalty1}%
  \advance\SplitTest-\insertpenalties}
\def\midinsert{\par
 \SaveLastSkip\penalty-150\SetSplitTest\RestoreLastSkip
 \ifnum\SplitTest=-1
  \@midfalse\p@gefalse\else\@midtrue\fi\@ins}
\def\@ins{\par\begingroup\setbox\z@\vbox\bgroup%
  \vglue\inserthardskipamount}
\def\endinsert{\egroup % finish the \vbox
  \if@mid \dimen@\ht\z@ \advance\dimen@\dp\z@
    \advance\dimen@\insertskipamount%            was 12pt (wn)
    \advance\dimen@\pagetotal\advance\dimen@-\pageshrink
    \ifdim\dimen@>\pagegoal\@midfalse\p@gefalse\fi\fi
  \if@mid%
    \ifdim\lastskip<\insertskipamount\removelastskip\insertskip\fi
    \nointerlineskip\box\z@\penalty-200\insertskip
  \else%
    \SaveLastSkip%                                  added (wn)
    \insert\topins{\penalty100 % floating insertion
    \splittopskip\z@skip
    \splitmaxdepth\maxdimen \floatingpenalty\z@
    \ifp@ge \dimen@\dp\z@
    \vbox to\vsize{\unvbox\z@\kern-\dimen@}% depth is zero
    \else \box\z@\nobreak\insertskip\fi}% was \bigskip\fi (wn)
    \RestoreLastSkip%                               added (wn)
   \fi\endgroup}
%% End Walter's insert stuff

 %%%%% Footnotes %%%%%

  \newcount\notenumber
  
  \def\note{\advance\notenumber by 1
    \footnote{\the\notenumber)}}

  \newbox\footbox

 %% The following modifies Plain TeX definitions, qv
  \def\footnote#1{\let\@sf\empty
    %{(the text)} is read later
    \ifhmode\edef\@sf{\spacefactor\the\spacefactor}\/\fi
    \sam${}^{\fam0 #1}$\@sf\vfootnote{#1}}%

  \def\vfootnote#1{\insert\footins\bgroup
     \interlinepenalty100 \splittopskip=1pt
     \floatingpenalty=20000
     \leftskip=0pt\rightskip=0pt%
     \parindent=.3em%% adjust
     \Smallfonts\rm%%osudeG added \Smallfonts
     \FootItem@{#1}%\strut% not nec
     \futurelet\next\fo@t}

  \def\FootItem@#1{\par\hangafter1\hangindent=\FootHang
     \setbox0=\hbox{\ignorespaces#1\unskip}%
     \dimen0=.4em\SetOverhang@% dimen0 is extra space
     \noindent\rlap{\box0}\kern\Overhang\ignorespaces}

  %\MaxFootTag{2)}%% in param file

  \def\fo@t{\ifcat\bgroup\noexpand\next \let\next\f@@t
    \else\let\next\f@t\fi \next}
  \def\f@@t{\bgroup\aftergroup\@foot\let\next}
  \def\f@t#1{\baselineskip=10pt\lineskip=1pt
            \lineskiplimit=0pt #1\@foot}%
     %%osudeG added \baselineskip=? pt\lineskiplimit=0pt
  \def\@foot{%%% special strut osu for end of each note
        \hbox{\vrule height0pt depth5pt width0pt}
        \egroup}
  \skip\footins=12 pt plus 0pt minus 0pt %% was \bigskipamount
    %% space added when footnote is present
  \count\footins=1000 % footnote magnification factor (1 to 1)
  \dimen\footins=8in % maximum footnotes per page

 %%%% Altenatives

  %%  Editorial stuff (delete??)

 \def\osumess#1{\EdSpider{\immediate\write16{Line \the\inputlineno: #1}}}%
 \def\HideEdStuff{\gdef\EdSpider##1{}}

 \font\BigSym=cmmi10 scaled \magstep 4

 \def\change{\InLMargin{\hbox{\BigSym \char63\kern10pt}}}

 \def\beginchange{\InLMargin{\hbox{\sam\twelvepoint$\heartsuit$\kern10pt}}}

 \def\endchange{\InLMargin{\hbox{\sam\twelvepoint$\spadesuit$\kern10pt}}}

 \def\InLMargin#1{\strut\vadjust{%
     \kern-\strutdepth
     \vtop to \strutdepth{%
         \baselineskip\strutdepth
         \llap{\sam$\smash{\hbox{\EdSpider{#1}}}$}\null}}}

 \def\strutdepth{\dp\strutbox}
 \def\strutheight{\ht\strutbox}

 \def\NoteInRMargin#1{\strut\vadjust{%
     \kern-1.001\strutdepth
     \vtop to \strutdepth{%
       \baselineskip\strutdepth
       \vss\rlap{\ninepoint\unskip\hskip\hsize
         \vtop to 0pt{%
           \hsize=16em\hfuzz=\hsize
           \leftskip=10pt%
           \rightskip=0pt plus 10000pt%
           \baselineskip=9.8pt\lineskip=.2pt%
           \let\\\break
           \noindent\EdSpider{#1}\vss}%
                \kern10pt}\hbox{}}%%\hbox{}=\null crucial!!
       }}

 \def\ednote#1{\NoteInRMargin{\tentt #1}}

 \def\cbar{\InLMargin{%
      \dimen0=\strutdepth\advance\dimen0 by \lineskip
      \vrule width 3pt
      height \strutheight depth \dimen0 \kern
      3pt}}

 \def\ccbar{\InLMargin{%
      \dimen0=2\strutdepth\advance\dimen0 by 2\lineskip
      \vrule width 3pt
        height 3\strutheight depth \dimen0 \kern
      3pt}}

 \newinsert\TRMargIns
 \dimen\TRMargIns=\maxdimen
 %\count\TRMargIns=0
 %\skip\TRMargIns=0pt

  \def\Ednote#1{\insert\TRMargIns{%
       \vbox to 0pt{\hsize=140pt\hfuzz=\hsize
           \leftskip=6pt%
           \rightskip=0pt plus 10000pt%
           \baselineskip=9.8pt\lineskip=.2pt%
           \let\\\break
           %\vglue\pagetotal% misplaces notes if inserts are present
           \SetPageRemainder% This ...
           \vglue540pt\vglue-\PageRemainder%  .. is a fix (WN)
           \noindent\EdSpider{\tentt #1}\vss}%
       \smallskip}}

 \def\KillEdStuff{\def\ednote##1{}\def\Ednote##1{}%
      \let\change\relax\let\beginchange\relax\let\endchange\relax
       \let\cbar\relax\let\ccbar\relax}

 %%% Compatibility with osumrip.sty
  %%

 %%% Parameters
  \topskip=12pt
  \newskip\StdBaselineskip % to set \baselineskip
  \StdBaselineskip 12pt
  \lineskip=1.1pt
  \lineskiplimit=.8pt
  \widowpenalty=10000 % 8000 to 10000
  \clubpenalty=10000  % 8000 to 10000
  \abovedisplayskip=6pt plus 1pt minus 1pt
  \abovedisplayshortskip=3pt plus 1.5pt
  \belowdisplayskip=6pt plus 1pt minus 1pt
  \belowdisplayshortskip=5pt plus 1pt minus 1pt
  \hfuzz=1.5pt   % Enable overfull box warnings at console

  \def\StdPretolerance{100}
  \tolerance=\StdPretolerance

  \newdimen\StdMathsurround
  \StdMathsurround=1.5pt % 1pt usual without \Mas
  \mathsurround=\StdMathsurround
  \Mas                   %% sophisticated mathsurround on
 % \Sam                   %% sophisticated mathsurround off

%% marker before English punctuation in displayed math
   \def\prose{\relax\hbox{\kern.6\StdMathsurround}}
  
  \def\StdParskip{0pt}    %% Larry wants {2pt plus 1pt}
  \parskip=\StdParskip
  \parindent=0.5cm
 
%%%% load Times for main body font

  \def\Times{ptmr  } 
  \def\TimesI{ptmri  } 
  \def\TimesB{ptmb  }
  \def\TimesBI{ptmbi  }
  \def\HelveticaN{phvrrn }

  =\Times at 10bp% roman text
  =\TimesB at 10bp% boldface extended
   % slanted roman
  \font\tenit=\TimesI at 10bp% text italic
  =\TimesBI at 10bp

  \font\tenmrm=cmr10  %%new name for math role at full size

%%%%% Fonts at ninepoint %%%%%

    =\Times at 9bp 
    \font\nineit=\TimesI at 9bp 
    =\TimesB at 9bp 
    =\TimesBI at 9bp 

    =\HelveticaN at 9bp 
       % see below

%%%%% Fonts at twelvepoint %%%%%

  =\Times at 12bp
  \font\twelveit=\TimesI at 12bp
  =\TimesB at 12bp

%%%%% Fonts at titlepoint %%%%%

  \font\titleit=\TimesI at 14.4bp
  =\TimesB at 14.4bp

 \SetAuthorHead{AuthorHead} % needs \ninepoint since box set
 \SetTitleHead{TitleHead}  % notably \HeaderFont

%%%% Char adjustments %%%%

  \def\lBr{\raise.125ex\hbox{[\kern.1125ex}}
  \def\rBr{\raise.125ex\hbox{\kern.1125ex]}}

 \setbox\footbox=\hbox{\Smallfonts 2)~}

%% Some optional font dimension and spacing 
%% adjustments beyond this point

%% Correct the lousy spacing of italic f (a hack).

  \bgroup
  \catcode`\@=11 %localised
  \gdef\itSpacing{%
     \xspaceskip=.31em plus.1em minus.05em \sfcode `f=2001
     \itWarning@\let\itWarning@\itWarning@@}
  \gdef\itSpacingOff{%
     \xspaceskip=0pt \sfcode `f=1000
     \let\itWarning@\relax}
   \global\let\itWarning@\relax
  \gdef\itWarning@@{\errmessage{%
  Special italic spacing already in force
  (you have probably omitted an ``endth'').
  See itSpacing macro in osuPSfnt.sty
         }}
  \egroup

 %%% Provisional fontdimen settings
  %%
 \fontdimen1\titlebf=0.0pt
 \fontdimen2\titlebf=3.6135pt
 \fontdimen3\titlebf=2.8908pt
 \fontdimen4\titlebf=1.44539pt
 \fontdimen5\titlebf=6.64882pt
 \fontdimen6\titlebf=14.45398pt
 \fontdimen7\titlebf=1.60439pt

 \fontdimen1\tenbi=0.26794pt
 \fontdimen2\tenbi=2.50937pt
 \fontdimen3\tenbi=2.00749pt
 \fontdimen4\tenbi=1.00374pt
 \fontdimen5\tenbi=4.59717pt
 \fontdimen6\tenbi=10.03749pt
 \fontdimen7\tenbi=1.11415pt

 \fontdimen1\twelverm=0.0pt
 \fontdimen2\twelverm=3.01125pt
 \fontdimen3\twelverm=2.409pt
 \fontdimen4\twelverm=1.2045pt
 \fontdimen5\twelverm=5.39615pt
 \fontdimen6\twelverm=12.045pt
 \fontdimen7\twelverm=1.33699pt

 \fontdimen1\twelveit=0.27731pt
 \fontdimen2\twelveit=3.01125pt
 \fontdimen3\twelveit=2.409pt
 \fontdimen4\twelveit=1.2045pt
 \fontdimen5\twelveit=5.37207pt
 \fontdimen6\twelveit=12.045pt
 \fontdimen7\twelveit=1.33699pt

 \fontdimen1\twelvebf=0.0pt
 \fontdimen2\twelvebf=3.01125pt
 \fontdimen3\twelvebf=2.409pt
 \fontdimen4\twelvebf=1.2045pt
 \fontdimen5\twelvebf=5.5407pt
 \fontdimen6\twelvebf=12.045pt
 \fontdimen7\twelvebf=1.33699pt

 \fontdimen1\tenrm=0.0pt
 \fontdimen2\tenrm=2.50937pt
 \fontdimen3\tenrm=2.00749pt
 \fontdimen4\tenrm=1.00374pt
 \fontdimen5\tenrm=4.49678pt
 \fontdimen6\tenrm=10.03749pt
 \fontdimen7\tenrm=1.11415pt

 \fontdimen1\tenit=0.27731pt
 \fontdimen2\tenit=2.50937pt
 \fontdimen3\tenit=2.00749pt
 \fontdimen4\tenit=1.00374pt
 \fontdimen5\tenit=4.47672pt
 \fontdimen6\tenit=10.03749pt
 \fontdimen7\tenit=1.11415pt

 \fontdimen1\tenbf=0.0pt
 \fontdimen2\tenbf=2.50937pt
 \fontdimen3\tenbf=2.00749pt
 \fontdimen4\tenbf=1.00374pt
 \fontdimen5\tenbf=4.61723pt
 \fontdimen6\tenbf=10.03749pt
 \fontdimen7\tenbf=1.11415pt

 \fontdimen1\ninerm=0.0pt
 \fontdimen2\ninerm=2.25842pt
 \fontdimen3\ninerm=1.80673pt
 \fontdimen4\ninerm=0.90337pt
 \fontdimen5\ninerm=4.0471pt
 \fontdimen6\ninerm=9.03374pt
 \fontdimen7\ninerm=1.00273pt

 \fontdimen1\nineit=0.27731pt
 \fontdimen2\nineit=2.25842pt
 \fontdimen3\nineit=1.80673pt
 \fontdimen4\nineit=0.90337pt
 \fontdimen5\nineit=4.02904pt
 \fontdimen6\nineit=9.03374pt
 \fontdimen7\nineit=1.00273pt

 \fontdimen1\ninebf=0.0pt
 \fontdimen2\ninebf=2.25842pt
 \fontdimen3\ninebf=1.80673pt
 \fontdimen4\ninebf=0.90337pt
 \fontdimen5\ninebf=4.15552pt
 \fontdimen6\ninebf=9.03374pt
 \fontdimen7\ninebf=1.00273pt

 %%% \SetExtraSpaces \MaxSpaceFactor \SetSpaceFactors
  %%  See TeXbook, page 76.

 \newcount\MaxSpaceFactor
 \MaxSpaceFactor=3000 %% to reset later

 %%%%% Tag styles and (hang-) indents
 \def\ItemStyle{\rm}
 \def\NrStyle{\rm}
 \def\ItemItemStyle{\rm}

 %% Analog dimensioning, convenient for local modifications:
 \MaxItemTag{(iii)}
 \MaxItemItemTag{(iii)}
 \MaxNrTag{(2)}
 \MaxFootTag{2)}
 % \MaxReferenceTag{AaaAA} % for biblio
 \def\ReferenceHang{30pt}

 \catcode`\@=\active

%%%%% End of hack of Neumann-Siebenmann macros

\loadbold

=\Times  
=\Times scaled750
=\Times scaled650
\font\rms=\Times scaled 920 

=\TimesBI scaled 860
=\TimesI scaled 860

\textfont0=\rrm  
\scriptfont0=\erm 
\scriptscriptfont0=\srm

\def\Augment#1#2{%
    \toks0\expandafter{#1}\toks2{#2}%
    \edef#1{\the\toks0\the\toks2}}

 \font\twelverma=\Times  scaled 1200
 \font\tenrma=\Times  scaled 1000
 \font\ninerma=\Times scaled 920
 =\Times scaled 840
 \font\sevenrma=\Times scaled 760
 =\Times scaled 680
 \font\fiverma=\Times scaled 600

 \Augment\tenpoint{%
  \textfont0=\tenrma  \scriptfont0=\sevenrma  
  \scriptscriptfont0=\fiverma  }

 \Augment\ninepoint{%
  \textfont0=\ninerma  \scriptfont0=\sevenrma 
  \scriptscriptfont0=\fiverma}

 \Augment\twelvepoint{%
  \textfont0=\twelverma  \scriptfont0=\ninerma  
  \scriptscriptfont0=\sevenrma}

\mathsurround=1pt
\hsize=13.45truecm
\vsize=19.5truecm
\hoffset=1.25truecm
\voffset=2truecm
\advance\baselineskip by 2pt

\predefine\til{\~}
\def\~#1{\relax\ifmmode\widetilde{#1}\else\til{#1}\fi}

\redefine \le{\leqslant}
\redefine \ge{\geqslant}
\define \wt#1{\mathaccent"0365{#1}}
\define \wh#1{\mathaccent"0362{#1}}

\define \iss{\,\Mathaccent{\raise -.8 ex\hbox{$\widetilde{}$\kern.1em}}\rightarrow\,}

\define \bigcupp{{\overset\cdot\to\bigcup}}

\define \prlim{{\varprojlim}\vphantom{i}\,}
\define \inlim{{\varinjlim}\vphantom{i}\,}

\define \alg{\mathop{\fam0 alg}}

\def\bitem{\item{$\bullet$}}

\Mas
\HideEdStuff
\rm 
 
%%%% For GT headers and footers:

\def\issn{{\nineit ISSN 1464-8997 (on line) 1464-8989 (printed)}}

\def\gtp{{\nineit Published 10 December 2000: \ \copyright\ Geometry \& 
Topology Publications}}

\def\gtv3{{\nineit Geometry \& Topology Monographs, Volume 3 (2000) --
Invitation to higher local fields}}

%%%%% For section idents:

\def\lione
{{\rms Geometry \& Topology Monographs}}

\def \litwo{{\rms Volume 3: Invitation to higher local fields
}} 

\def\tinfo #1.#2.#3-#4
{{
\noindent  {\lione} \hfill 
\par 
\vskip-1.5pt
\noindent {\litwo} \hfill
\par 
\vskip-1,5pt
\noindent {\rms Part #1, section #2, pages #3--#4} \hfill
\vskip24pt 
}}

\def\tinfos #1.#2.#3-#4
{{
\noindent  {\lione} \hfill 
\par 
\vskip-1.5pt
\noindent {\litwo} \hfill
\par 
\vskip-1.5pt
\noindent {\rms Pages #3--#4} \hfill
\vskip24pt 
}}

\def\tinfoi #1
{{
\noindent  {\lione} \hfill 
\par 
\vskip-1.5pt
\noindent {\litwo} \hfill
\par 
\vskip-1.5pt
\noindent {\rms Pages iii--xi: Introduction and contents} \hfill
\vskip26pt 
}}

%%%% Set headers and footers %%%%

  \def\titlepagehead{\hfil}

  \newif\iftitlepage\titlepagefalse
  \newif\ifblankpage\blankpagefalse
  \def\makeheadline{
     \ifblankpage{}\else%
     \iftitlepage
\vbox{\line{\vbox to 8.5pt{}
\ninerm
\copy\HLinebox \hfill
\hglue5mm\ninebf\folio 
\titlepagehead}}%
      \else
\vbox{\ifodd\pageno\rightheadline\else\leftheadline\fi}%
      \fi\vskip 12pt\fi}%
     \def\rightheadline{\line{\vbox to 8.5pt{}%
      \ninerm
\copy\TitleBox \hfill
\hglue5mm\ninebf\folio}}%
     \def\leftheadline{\line{\vbox to 8.5pt{}%
        \unskip\ninerm\unskip\ninebf\folio\hglue5mm
      %*%
 \hfill \copy\AuthorBox
%\hfill
}}

 \footline={\ifblankpage{}\else
\iftitlepage\ninepoint\sam\hfill%} 
\line{\vbox to 8.5pt{}%\ninerm
\copy\TFLinebox
\hfill
\hglue5mm %\ninebf\folio
}
            \else
\ninepoint\sam\hfill%}
\line{\vbox to 8.5pt{}%\ninerm
\copy\FLinebox
\hfill 
\hglue5mm
}
\hfil\fi\global\titlepagefalse\fi}

\def\blankpage{{\blankpagetrue\noindent\hbox to 10pt{\hss}\vfill
\pagebreak}}

\tenpoint\rm %% always start here
 
  %%% all done and macros loaded!